\documentclass[brochure,english,12pt]{smfbourbaki}

\usepackage{lmodern,amssymb,bm}

\usepackage{babel}
\usepackage{graphicx}
\usepackage{mathtools}
\usepackage[colorlinks=true, linkcolor=blue, citecolor=red, urlcolor=blue]{hyperref}

\usepackage{cancel}

\usepackage[stretch=10,shrink=10,step=2,kerning=true,protrusion=true,expansion=true,final]{microtype} 

\usepackage[
backend=bibtex,
style=authoryear, 
citestyle=authoryear-comp,
maxnames=7,
sortcites=false %%%% to keep the order in \cite command
]{biblatex}
\usepackage{csquotes}

\newtoggle{tnbcbx@howcited}
\DeclareEntryOption[boolean]{howcited}[true]{\settoggle{tnbcbx@howcited}{#1}}
\DeclareBibliographyOption[boolean]{howcited}[true]{\settoggle{tnbcbx@howcited}{#1}}
\DeclareTypeOption[boolean]{howcited}[true]{\settoggle{tnbcbx@howcited}{#1}}

\newbibmacro{howcited}{%
  \iftoggle{tnbcbx@howcited}
    {\iffieldundef{shorthand}
       {}
       {\setunit{\addspace}%
        \printtext[parens]{%
          \bibstring{citedas}%
          \setunit{\addcolon\space}%
          \printfield{shorthand}%
          \setunit{\addslash}%
          }}}
    {}}

\renewbibmacro{finentry}{\usebibmacro{howcited}\finentry}

\DefineBibliographyStrings{french}{citedas    = {cit\'e ci-dessus avec l'acronyme}}
\DefineBibliographyStrings{english}{citedas    = {heretofore cited as}}

\DefineBibliographyExtras{french}{\restorecommand\mkbibnamefamily}

\DeclareDelimFormat{nameyeardelim}{\addcomma\space}

\DeclareNameAlias{sortname}{given-family}

\renewcommand{\bibnamedash}{\leavevmode\raise3pt\hbox to3em{\hrulefill}\space}

\AtEveryBibitem{%
  \clearfield{issn} % Remove issn
  \clearfield{isbn} % Remove isbn
  \clearfield{doi} % Remove doi
  \clearlist{language} %%% remove language
  \ifentrytype{online}{}{% Remove url except for @online
  \ifentrytype{unpublished}{}{% ou bien @unpublished
    \clearfield{url}
  }
  }
}

\renewbibmacro{in:}{%
    \ifentrytype{article}{}{\printtext{\bibstring{in}\intitlepunct}}}

\DeclareFieldFormat[article,periodical,inreference]{number}{\mkbibparens{#1}}
\DeclareFieldFormat[article,periodical,inreference]{volume}{\mkbibbold{#1}}
\renewbibmacro*{volume+number+eid}{%
    \printfield{volume}%
    \setunit*{\addthinspace}% NEW
    \printfield{number}%
    \setunit{\addcomma\space}%
    \printfield{eid}}

\DeclareFieldFormat[article,inbook,incollection]{title}{\enquote{#1}\addcomma} 

\date{Janvier 2026}
\bbkannee{78\textsuperscript{e} ann\'ee, 2025--2026}  % Commandes sp\'eciales
\bbknumero{1247}                                      % Commandes sp\'eciales

\newcommand{\F}{f^\varepsilon}
\newcommand{\Z}{\mathbf{z}}
\newcommand{\X}{\mathbf{x}}
\newcommand{\V}{\mathbf{v}}
\newcommand{\bbM}{\mathbb{M}}

 \def\eps {{\varepsilon}}
\def\d {{\partial}}
\def\indc {{\textbf 1}}

\title{Derivation of the  Boltzmann equation \\
from hard-sphere dynamics}
\subtitle{after Y. Deng, Z. Hani, and X. Ma}

\author{Thierry Bodineau}
\address{I.H.E.S., Universit\'e Paris-Saclay,  CNRS \\  Laboratoire Alexandre Grothendieck,\\
35 Route de Chartres \\91440 Bures-sur-Yvette, France}% Email address
\email{bodineau@ihes.fr}

\author{Isabelle Gallagher}
\address{Université Paris Cité, Sorbonne Université \\ CNRS, IMJ-PRG \\75013 Paris, France}
\email{isabelle.gallagher@u-pariscite.fr}

\author{Laure Saint-Raymond}
\address{I.H.E.S., Universit\'e Paris-Saclay,  CNRS \\  Laboratoire Alexandre Grothendieck,\\
35 Route de Chartres \\91440 Bures-sur-Yvette, France}
\email{laure@ihes.fr}

\author{Sergio Simonella}
\address{Sapienza Universit\`a di Roma \\Dipartimento di Matematica G. Castelnuovo \\ Piazzale A. Moro 5, 00185 Roma, Italy}
\email{sergio.simonella@uniroma1.it}

\begin{document}

\maketitle

\tableofcontents

\section*{Introduction}

The Boltzmann equation describes the time evolution of the  particle density  of a collisional rarefied gas in the position and velocity phase space. Its success, together with the kinetic theory it gave rise to, represents an undisputed achievement in physics.

In this note we present the result by  \textcite{DHM}, who derive the Boltzmann equation  starting from a microscopic system of hard spheres, for  times much larger than the mean free time (the average time between two collisions of a given particle): the derivation holds  on the same time interval on which there is a smooth solution of the Boltzmann equation.

Let us start by describing the main theorem informally. A precise statement requires a number of definitions  and will appear on page \pageref{thm: DHM} (Theorem \ref{thm: DHM}). Consider~$N$ hard spheres moving in $\mathbb R^d$, initially ``independent" and identically distributed according to a distribution $f^0= f^0(x,v)$. Assume that the Boltzmann equation with initial data $f^0$ admits a smooth solution $f$ on a time-span $[0,T]$. Then, for any $t \in [0,T]$, the (random) fraction of hard spheres in any set $\Delta \subset \mathbb R^{2d}$ converges in probability to~$\displaystyle \int_{\Delta} f(t)$, in the limit when the number of particles $N \to \infty$ and their diameter $\eps \to 0$ in such a way that the mean free time is of order 1.

{
The problem of deriving the Boltzmann equation from microscopic systems has been addressed by  many authors in the past century.
%, starting with the pioneering works  of \textcite{Bogoliubov}, \textcite{BG},  \textcite{Kirkwood}    and \textcite{Yvon}, who gave their names to the BBGKY hierarchy of equations (for smooth potentials) on the successive collection  of correlation functions in the 1940s\footnote{{   In statistical mechanics, correlation functions are now a customary way to describe the joint probability distribution of finitely many particles in phase space.}}. This was followed by 
An initial key step is}
 the enlightening work of \textcite{Grad} on the low-density scaling, where he 
 % derived {  heuristically} 
 {wrote down}
an effective differential equation for the first correlation function in the low-density regime. 
 The first mathematical result on this problem goes back to \textcite{Cercignani}  and to the celebrated theorem of~\textcite{Lanford},  who proved that in the special case of the hard-sphere system, the asymptotic factorisation of the correlation functions (known in kinetic theory as the {\it propagation of chaos} property)  can be established by a careful study of trajectories. In particular,      Lanford  demonstrated that there is no contradiction between the reversibility in the microscopic
laws of mechanics and the irreversibility exhibited by Boltzmann's equation. The proof  was an important breakthrough and after a number of subsequent works until the early 2000's,  the details of the argument were completed {and the result refined and generalized\footnote{See~\textcite{Spohn91},~\textcite{CIP94},~\textcite{CGP},~\textcite{MT12}~\textcite{GSRT14},~\textcite{PSS14},~\textcite{Ayi},~\textcite{GP21},~\textcite{LeBihan22},~\textcite{Dolmaire},~\textcite{DolmaireSaffirio}.}.}  The main limitation of the method  ---  unsurpassed until the work of  Deng, Hani and Ma except for some specific settings\footnote{{Perturbation of vacuum in \parencite{IP, Denlinger}, tagged particle at equilibrium in (van \cite{vBLLS}; \cite{BGSR}; \cite{Catapano, Fougeres}) or fluctuations around equilibrium in \parencite{BGSRS23-1, BGSRS24, LeBihan25, LeBihan25-1}}.}  is that the result is only valid for a very short time:  essentially a  {small} fraction of the mean free time\footnote{ {To keep in mind the order of  physical magnitudes, for a gas at room temperature and atmospheric pressure
the mean free time is of the order of $10^{-9}$ seconds, while the mean free path is of the order of $10^{-5}$ meters.}}.
The reason for this drawback  is that the method of proof consists in expanding the correlation functions into    time series, which only converge for short times. 
Actually, the Boltzmann equation is itself very difficult to study mathematically and the existence and uniqueness of solutions for long times is only known in specific settings. This lack of  stability is due to the locality of the interaction, which is also one of the reasons why it is so difficult to prove the convergence of the microscopic dynamics.

Deng, Hani and Ma were able, fifty years after the seminal paper by Lanford,  to go beyond  the short time restriction and to prove the long-awaited result that the convergence of the hard-sphere dynamics holds up to any time for which  the Boltzmann equation has a smooth solution.
Important ideas for this proof have already been devised in the series of works by
\textcite  {DHw1,DHw2,DH-waves} to derive a kinetic equation in the context of wave turbulence on long times
(see Remark \ref{rem: wave} for some comments).

In this paper, we present some of the key ideas of the proof of \textcite{DHM}. In Section~\ref{sec:Boltzmann}
we introduce the Boltzmann equation and  present some classical well-posedness results and their limitations. Section~\ref{sec:STD} is devoted to the presentation of the short time derivation result by  {Lanford} (Theorem~\ref{thm: Lanford}). 
 A pictorial sketch of proof  of
Lanford's theorem is provided in Section~\ref{sec:proofL}: this  somewhat unusual presentation   is designed as a preparation to the presentation of the proof of the  theorem   by Deng, Hani and Ma (Theorem~\ref{thm: DHM}), which is provided in Section~\ref{sct:DHM}. Some consequences of the result, and perspectives, are given in Section~\ref{sct:perspectives}.

 \paragraph*{Acknowledgements}
 We are very grateful to Fran{\c c}ois Golse, Mario Pulvirenti, Herbert Spohn and Isabelle Tristani  for useful comments on a previous version of this manuscript.  We especially thank  Nicolas Bourbaki for a very thorough reading and many suggestions. S.S. acknowledges funding from the European Union (ERC CoG KiLiM, 101125162).

\section{The Boltzmann equation}
\label{sec:Boltzmann}

The Boltzmann equation,   devised by  {Boltzmann} in 1872, is the cornerstone of kinetic theory. It provides a statistical description of  a gas constituted of    identical particles moving in the~$d$-dimensional space~$\mathbb R^d$ ($d=2,3$)
with uniform rectilinear motion between each elastic binary collision (see Section~\ref{sct:hard spheres}
 below for more).

Under the assumption that the particles are independent and identically distributed (this is referred to as the ``chaos assumption"), their {common  distribution $f$  is expected to be a real-valued function
depending on $t \in \mathbb R$ (time), $x \in \mathbb R^d$ (position) and $v \in \mathbb R^d$ (velocity) satisfying the Boltzmann equation}
\begin{equation}
\left\{ \begin{aligned}
&\underbrace{\partial_t f(t,x,v) \! +\!v \cdot \nabla _x f (t,x,v)}_{\tiny\mbox{transport}}
\!= \! \displaystyle\int_{\mathbb R^{d} \times\mathbb S^{d-1} }   \! \! \big( \underbrace{f(t,x,v_\star') f(t,x,v')}_{\tiny\mbox{gain term}} - \underbrace{f(t,x,v_\star) f(t,x,v)}_{\tiny\mbox{loss term}}\big) 
\\
&  \qquad\qquad \quad\qquad \qquad\qquad\qquad\qquad\qquad\qquad  \qquad \times \big ((v-v_\star)\cdot \omega\big)_+ \, dv_\star\, d\omega \,\,  ,\\
&  f(0,x,v) = f^0(x,v)
 \end{aligned}
  \right. 
\label{eq:Beq}
  \end{equation}
where the ``pre-collisional'' velocities $(v',v_\star ')$ in the ``gain'' term on the right-hand side are defined by the scattering laws
\begin{equation}
\label{eq: scattlaw}
 v' \coloneqq v- \big( (v-v_\star) \cdot \omega\big)\,  \omega \,  ,\qquad
 v_\star' \coloneqq v_\star+\big((v-v_\star) \cdot \omega\big)\,  \omega  \, .
\end{equation}
 \begin{figure}[h] 
\centering
\includegraphics[width=1.5in]{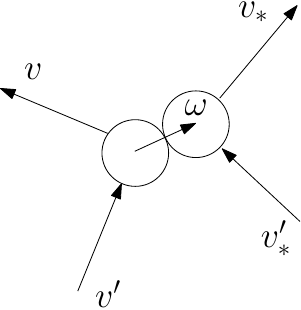} 
\caption{\small 
The hard-sphere scattering law.}
\label{HS-scatt}
\end{figure} 

Note that
\begin{equation}\label{conservation momentum energy}
	\begin{aligned}
		v+v_\star & = v'+v_\star'  \,,\\
		|v|^2+|v_\star|^2 & = |v'|^2+|v_\star'|^2 \,,
	\end{aligned}
\end{equation}
which express the fact that  collisions are  elastic   and thus momentum and energy are
conserved during collisions. The initial value of~$f$ is prescribed to be a given function~$f^0$. We have used the notation~$a_+$ for~$\max(0,a)$.

In the following we shall use the shorthand notation  
$$f_\star \coloneqq  f(v_\star) \, , \quad f'\coloneqq f(v')  \quad \mbox{and} \quad f'_\star\coloneqq f(v_\star') \, ,
$$
so that the right-hand side of~(\ref{eq:Beq})   becomes
\begin{equation} \label{eq:BCO}
Q(f,f)(v)\coloneqq\int \big( f'f_\star'-ff_\star\big) \big ((v-v_\star )\cdot \omega\big)_+ \,  dv_\star\,  d\omega\,.
\end{equation}

\subsection{Conservation laws and H-Theorem}  
From~\eqref{conservation momentum energy} and using the well-known fact (see for instance~\cite{Cercignani}) that the mappings~$(v,v_\star)\mapsto(v_\star,v)$ and~$(v,v_\star,\omega)\mapsto \left(v',v_\star',\omega\right)$  have  unit Jacobian determinants, one can show   formally that for any test function~$\varphi$ defined on~$\mathbb R^d$, and with the notation introduced at the end of the previous paragraph, 
$$
	\int Q(f,f)\,\varphi\, dv=\frac14 \int [f'f'_\star-ff_\star] (\varphi+\varphi_\star - \varphi'-\varphi'_\star) \,\big ((v-v_\star)\cdot \omega\big)_+   \, dv\,dv_\star\,d\omega\,.
$$
In particular, it can be shown that
$$	\int Q(f,f)\varphi\,dv=0 \quad \mbox{for all $f$ regular enough} 
$$
 if and only if $\varphi(v)$ is a collision invariant, meaning that   $\varphi(v)$ is a linear combination of~$\{1,v^1,\dots, v^d,|v|^2\}$. Thus, successively multiplying the Boltzmann equation~(\ref{eq:Beq}) by the collision invariants and then integrating in  velocity gives rise formally to  the {\it local conservation laws}:
$$
\begin{aligned}
\d_t \int  f(t,x,v) \,dv+ \nabla_x \cdot  \int  f(t,x,v) v\,dv = 0& \,;\\
\d_t \int  f(t,x,v) v\,dv+ \nabla_x \cdot  \int  f(t,x,v) v \otimes v \,dv = 0 & \,;\\
\d_t \int  f(t,x,v) \frac{|v|^2}2 \,dv+ \nabla_x \cdot  \int  f(t,x,v) \frac{|v|^2}2 v \,dv = 0 & \, .
 \end{aligned}
$$
  The second line should be understood component-wise as
$$
  \forall 1 \leq i \leq d \, , \quad \d_t \int  f(t,x,v) v^i\,dv+  \sum_{j = 1}^d \partial_j \int  f(t,x,v) v^iv^j\,dv = 0\, .
$$
Taking $\varphi = \log f$ and using the additivity of the logarithm shows that the collision operator  is {\it dissipative} in the sense that
$$
D(f) \coloneqq  -\int Q(f,f) \log f  \,dx \,dv\geq 0 $$
with equality if and only if~$\log f$ is a  collision invariant. 
This  leads to the entropy inequality
\begin{equation} \int f\log f (t,x,v) \,dx\,dv +\int_0^t D(f)(s)\,ds \leq \int f^0 \log f^0 (x,v) \,dx\,dv\,,
 \label{eq:Hthm}
 \end{equation}
known since Boltzmann as the H-Theorem. Note that this increase in the entropy, defined by~$\displaystyle\int - f\log f (t,x,v) dxdv$ (which is strict unless~$f$ is an equilibrium, i.e.\,a spatially homogeneous Gaussian distribution) reflects the fact that Boltzmann's equation is   {\it time  irreversible}, in contrast with the particle system: we refer to Paragraph~\ref{Loschmidt} below for more on this point.

%In particular due to the above considerations, $D(f)$  vanishes identically if and only if~$\log f$ is a linear combination of collision invariants. Let us denote  the   Maxwellian 
%$$
% M_\beta (v) \coloneqq \frac{\beta^{d/2}}{(2\pi)^\frac d2} \exp \left( - \beta \frac{ |v|^2}{2} \right)  \,, 
% $$
%which is  an equilibrium in the sense that it cancels each term of~(\ref{eq:Beq}).
%Then, still in a formal sense, 
%\begin{equation}\label{H-Theorem}
%\begin{aligned}
% H(f|M_\beta) (t) &\coloneqq \int    \left( f \log {f\over M_\beta} -f+M_\beta \right) (t,x,v)
%    dv  dx\\
%    & = H(f^0|M_\beta) - D(f)(t)\\
%    & \leq H(f^0|M_\beta) 
%\end{aligned}
%\end{equation}
%which means that~{\it the relative entropy decreases}. This is 

\subsection{Existence theories}  \label{sct:existence}
 In this paragraph we give a short presentation of three approaches to solving the Cauchy problem~(\ref{eq:Beq}): the question is to know whether given~$f^0$,  there exists  a unique  solution to the equation, and what is its life span.

\subsubsection*{Renormalised solutions}
For many partial differential evolution equations, {a natural approach to constructing solutions is to exploit the a priori controls provided by the equation in order to design an approximation scheme}. 

 We have presented in the previous paragraph the conservation laws associated with the Boltzmann equation~(\ref{eq:Beq}).  Unfortunately
such   a priori information  on the solution, typically the log-linear bound~(\ref{eq:Hthm}) ---  also known as an~$L\log L$ bound ---  is not sufficient to make sense of the collision term  which is quadratic. By introducing  a notion of {\it renormalised solution} (a function with bounded entropy satisfying a formally equivalent equation obtained by taming concentrations) and exploiting   the entropy inequality~(\ref{eq:Hthm}), \textcite{DiPL} nevertheless succeeded  in proving the {\it  global existence} of  renormalised solutions to~(\ref{eq:Beq}). However, their uniqueness   is not known to date, nor the fact that they satisfy (\ref{eq:Beq}) even in the sense of distributions.

\subsubsection*{Fixed-point methods}
 In such a situation where the a priori bound does not provide enough information to conclude, the strategy usually employed  is to change perspective and consider the PDE as a fixed point problem to be solved. More precisely,  using  {the so-called} Duhamel formula we can rewrite the equation~(\ref{eq:Beq}) {in integrated form}:
 \begin{equation}\label{eq:BeqFP}
\begin{aligned}
f(t) &= \mathcal S_t f^0 + \mathcal B_t\big (f(t),f(t)\big) , \\
  \mathcal B_t\big (f(t),f(t)\big) &\coloneqq\int_0^t \mathcal S_{t-s} Q \big (f(s), f(s)\big) \, ds,
\end{aligned}
\end{equation}
where $\mathcal S_t$ is the free transport operator  that attaches  to any function~$g$ on~$\mathbb R^{2d}$ the function
$$
\mathcal S_t g(x,v) \coloneqq g(x-vt,v) \, .
$$
Solving~(\ref{eq:BeqFP})   reduces to  looking for a Banach space in which the operators~$\mathcal S_t$ and $\mathcal B_t$ operate continuously. 
A typical result, based on the Cauchy-Kowalewskaya theorem (as in~\textcite{Treves} or~\textcite{Nirenberg} for instance) states that if $f^0$ is a continuous function in a weighted space~$L^\infty_\beta$ in velocities, meaning  that~$f^0$ is continuous on~$\mathbb R^{2d}$ and satisfies
\begin{equation}\label{eq:weighted space}
\| f_0\|_{L^\infty_\beta}\coloneqq \Big\|f_0\exp( \frac \beta 2 |v|^2)\Big \|_{L^\infty(\mathbb R^{2d})} <\infty
\end{equation}
for  some~$\beta>0$,  then there exists a constant $C_\beta>0$ (depending only on $\beta$) such that the solution of (\ref{eq:BeqFP})  has a unique continuous solution on~$[0,T]$, with
$$
T = \frac{ C_\beta }{\ \ \, \Big \|f_0\exp( \frac \beta 2 |v|^2)\Big\|_{L^\infty(\mathbb R^{2d})} } \, \cdotp
$$
 Actually not only is the solution continuous, but it is also stable with respect to perturbations of the initial data, in the same function space. We refer for instance to~\textcite{KanielShinbrot} for a result of this type. 
Unlike the renormalised solutions framework presented above, this type of result does not make any use of the particular structure of the nonlinear term (especially of the cancellations between the gain and loss terms): it suffices for the nonlinear term to be quadratic and to produce   a linear loss in velocity 
%{   -- actually  the use of the Cauchy-Kowalewskaya theorem is due to this loss in velocities}. 
 (actually this loss is triggered by the defect of integrability in~\eqref{eq:BCO} for large relative velocities $v-v_*$).

To exaggerate somewhat, by this method one is roughly   solving an equation of the type
$$
\partial_t f = f^2 \, .
$$ 
Such an equation is of course known to have a solution blowing up in finite time, so the only way to improve the well-posedness of (\ref{eq:BeqFP}) is to use more information on the structure of the nonlinear term.

\subsubsection*{Perturbative methods}
Using the entropy seems   difficult  as it provides too little a priori smoothness and decay on the solution, so a third strategy consists in   exploiting known global solutions and then to  act perturbatively. For instance, since the Maxwellian
$$
 M_\beta (v) \coloneqq \frac{\beta^{d/2}}{(2\pi)^\frac d2} \exp \left( - \beta \frac{ |v|^2}{2} \right)  
 $$
 is a stationary solution, one can   try to solve  (\ref{eq:BeqFP})  in the vicinity of~$M_\beta$. Assuming  the initial data takes the form
 $$
 f^0 = M_\beta(1+g^0)
 $$
 with~$g^0$ smooth and small enough in an appropriate function space (with strong enough decay in velocities), then one can  produce a unique solution~$f(t) =  M_\beta(1+g(t))
$ {\it for all positive times} \parencite{Ukai}. 
This approach crucially uses the fact that the collision operator $2 Q(M_\beta, \cdot )$ is coercive in $L^2 (M_\beta^{-1} dv)$, which is inherited from the dissipative nature of the equation. Another possibility is considering solutions which do not depend on the space variable \parencite{Carleman}.

\subsection{The Loschmidt paradox and the loss of information}  \label{Loschmidt}

What went down in history as the most famous criticism of Boltzmann's theory is an argument due to  \textcite{Loschmidt} (see \cite{Cercignani_book}, for an historical account). It confronts two indisputable facts: the microscopic laws of classical mechanics are reversible, yet the kinetic world displays an arrow of time, as the H-Theorem indeed shows. 

The apparent contrast can be interpreted as  related to  information: the irreversible evolution  corresponds to a gradual loss of microscopic detail. Entropy measures the level of information, and the entropy production  quantifies the amount of information discarded when the complete microscopic state of 
the particle system is replaced by a reduced statistical description.

In fact, a mathematical investigation of the transition from particle systems to the kinetic description allows one to complement this discussion with a more pragmatic perspective. It allows indeed to quantify the defect of chaos (see Section~\ref{sct:chaos} below), and to relate the entropy production to the growth of microscopic correlations (see \cite{Spohn97}).

\section{Short time derivation} \label{sec:STD}
In this section  we shall present Lanford's theorem, which ensures the validity of the Boltzmann equation~(\ref{eq:Beq}) starting with a system of hard spheres, for a short time. It also provides a  control on correlations, and therefore on the chaos approximation, i.e. the asymptotic factorisation of the correlation functions.

 \subsection{The hard-sphere system and the Liouville equation}  \label{sct:hard spheres}
\begin{figure}[h] 
\centering
\includegraphics[width=4in]{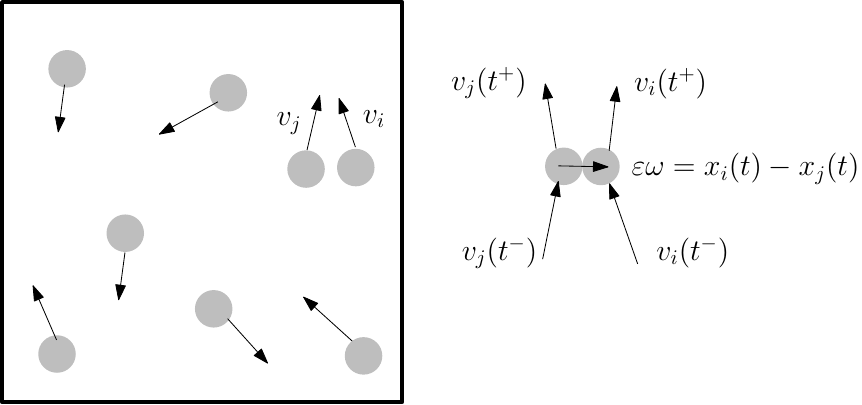} 
\caption{\small 
The hard-sphere dynamics and the deflection condition. }
\label{HS-fig}
\end{figure}

We
assume throughout that the gas is made of small identical particles, which are spheres of diameter~$\varepsilon \ll 1$
 interacting only through contact (see Figure \ref{HS-fig}). This is referred to as a \emph{hard-sphere gas}.    
If the gas consists of $N$ particles, its state is described by  the collection~$(x_1,\dots,x_N) \in \mathbb R^{dN}$ of  their positions and the collection~$(v_1,\dots,v_N) \in \mathbb R^{dN}$ of their velocities. 
Knowing  the state of the gas at a certain instant (say $t = 0$), we can  describe its state   at any time $t \in \mathbb R$ by solving the following system of ordinary differential equations:
\begin{equation}
\label{Newton} \frac {dx_i}  {dt} =  v_i\,,\quad  \frac {dv_i}  {dt} =0 \quad \hbox{as long as \ } |x_i(t)-x_j(t)|>\varepsilon 
\quad \hbox{for \ } 1 \leq i \neq j \leq N
\, ,
\end{equation}
with specular reflection condition at each collision: if~$|x_i(t)-x_j(t)|=\varepsilon$, then
\begin{equation}
\label{defZ'nij}
\begin{cases}
 v_i (t^+)\coloneqq v_i (t^-)-  \big( (v_i (t^-) - v_j (t^-))\cdot {\omega} \big) \, {\omega}   \\
 v_j(t^+) \coloneqq v_j(t^-)+ \big( (v_i (t^-) - v_j (t^-))\cdot {\omega} \big) \, {\omega} 
\end{cases}
\end{equation}
where ${\omega} \coloneqq (x_i (t) - x_j (t))/\varepsilon$ is the unit vector in the direction of the relative positions  at   the collision time~$t$. As pointed out in~(\ref{eq: scattlaw})-(\ref{conservation momentum energy}), this condition ensures that the collisions are elastic: momentum and energy are conserved. Note that~(\ref{Newton}) has a unique global solution~$\Big((x_1,\dots,x_N)(t),(v_1,\dots,v_N)(t)\Big)$ outside a set of initial data which has Lebesgue measure zero  (see~\cite{Alexander}, for example).

\bigskip
The dynamics of the particles is   deterministic, and we choose the initial state of the gas  randomly following a grand-canonical formalism (see, for example, \cite{Ruelle}). The grand-canonical formalism is known to be technically convenient when dealing with correlations at low density (see van~\cite{vBLLS}).

We therefore assume that the number of particles~$\mathcal N$ in the gas is random, with a Poisson-like distribution with parameter~$\mu_\eps>0$ {to be determined}.  Given~$N$ particles, we denote by~${\mathbf x}_N\coloneqq(x_1,\dots,x_N) $ in~$  \mathbb R^{dN}$   the vector representing the position of each of these particles and by~${\mathbf v}_N\coloneqq(v_1,\dots,v_N)  $ in~$  \mathbb R^{dN}$ the vector representing their velocities. We set~${\mathbf z}_N\coloneqq({\mathbf x}_N,{\mathbf v}_N)$.  As a shorthand notation we shall write also~$\Z_H$ for the set~$(z_i)_{i \in H}$ when~$H$ is a subset of~$\{1,\dots, N\}$.  For simplicity let us assume that the
initial density distribution~$W^{\varepsilon, 0}_N$ corresponds to almost independent and identically distributed particles
% of the particles is   then a function of $\Z_N\coloneqq(X_N,V_N)$, invariant under permutation
% since the particles are   identical.  The particles are not allowed to overlap, so~$W^{\varepsilon0}_N$  is  defined in a subspace of~$\mathbb R^{2dN}$ which we will denote
%$$
%\mathcal D^\varepsilon_N\coloneqq \Big\{\Z_N \in \mathbb R^{2dN} \, /\, |x_i-x_j|>\varepsilon \, \, \forall i \neq j\Big\} \, .
%$$
%Note that by the exclusion condition, $\mathcal D^\varepsilon_N$ is empty as soon as $N$ is sufficiently large.
%
%In this presentation, the probability density of having a configuration with $N$ particles at configuration~$\Z_N$ at time zero will be deliberately chosen as simple as possible (although more general results can be proved). Calling~$\mu_\varepsilon$ the scaling parameter, specified below in~(\ref{defmueps}):
\begin{equation}
\label{eq: initial measure}
\frac{1}{N!} W^{\varepsilon, 0 }_{N}(\Z_N) 
\coloneqq \frac{1}{\mathcal Z^\varepsilon} \,\frac{\mu_\varepsilon^N}{N!} \, 
{\textbf 1}_{\mathcal D^{\,\varepsilon}_{\!N} } (\Z_N) \; \prod_{i=1}^N f^0 (z_i) \, ,
\end{equation} 
where $\mathcal D^{\,\varepsilon}_{\!N} $ encodes the exclusion condition
$$
\mathcal D^{\,\varepsilon}_{\!N} \coloneqq \Big\{\Z_N \in \mathbb R^{2dN} \, \colon\, |x_i-x_j|>\varepsilon \, \, \forall i \neq j\Big\} $$
and the  partition function is defined by
\begin{equation}
\label{eq: partition function}
\mathcal Z^\varepsilon \coloneqq  1 + \sum_{N\geq 1}\frac{\mu_\varepsilon^N}{N!}  
\int_{\mathcal D^{\,\varepsilon}_{\!N} } d \Z_N \prod_{i=1}^N  f^0(z_i)  \,.
\end{equation}
Note that the small correlation due to the initial exclusion  is harmless in the analysis that will follow and we will ignore it most of the time.
The distribution of a single particle $f^0$ is assumed to be a continuous and Lipschitz probability density  on $\mathbb R^{2d}$ satisfying (\ref{eq:weighted space}) for some~$\beta>0$.

The probability of an event $A$ relative to the measure $\eqref{eq: initial measure}$ 
will be denoted by ${\mathbb P}_\varepsilon(A)$, and ${\mathbb E}_\varepsilon$ will be the expectation. In particular, the expectation of the number of particles~$\mathcal N$ is determined  asymptotically by $\mu_\varepsilon$ in the sense that $$
\lim_{\varepsilon \to 0}\mu_\varepsilon^{-1} {\mathbb E}_\varepsilon [\mathcal N] = 1 \, .
$$
We  assume from now on that the gas is {\it dilute} by setting
\begin{equation}
\label{defmueps}
 \mu_\varepsilon =\varepsilon^{-(d-1)} \, .
\end{equation}
Since the size of a cylinder based on a sphere of diameter~$\varepsilon$ is proportional to~$\varepsilon^{d-1} $,  this choice of scaling guarantees that the mean free time between collisions of a given particle is of the order of one, meaning that it is essentially independent of~$\varepsilon$. 

\bigskip
The Newtonian system satisfied by the particles translates, at the level of the distribution function~$W^{\varepsilon }_{N}(t)$ for all times~$t$, into the following Liouville equation 
\begin{equation}
\label{Liouville}
    \partial_t W^{\varepsilon}_N +\V_N \cdot \nabla_{\X_N} W^{\varepsilon}_N =0  \,\,\,\,\,\,\,\,\, \hbox{in } \,\,\,\mathcal D^{\,\varepsilon}_{\!N} \, ,
\end{equation}
	with specular reflection condition  at the boundary of~$\mathcal D^{\,\varepsilon}_{\!N} $.

	 \subsection{Lanford's theorem}   \label{sec:Lanford}
	 
We are interested in the typical behaviour of a particle in the hard-sphere gas. 
We recall that, according to the grand-canonical description given in Section \ref{sct:hard spheres}, the randomness comes through the initial configuration of particles $\left( z_1(0), \cdots, z_{\mathcal N}(0)\right)$, chosen according to the probability measure \eqref{eq: initial measure}. The main object of investigation becomes then the discrete random field
\begin{equation}
\label{eq: def empirical}
\Pi^{\varepsilon}(h,t) \coloneqq  \frac1{\mu_\varepsilon } \sum_{i=1}^\mathcal N h(z_i(t))
\end{equation}
defined on bounded test functions $h\in C^0(\mathbb R^{2d})$.
This random field is introduced, by following the original intuition of Boltzmann, to ``count'' the number of particles in little boxes of the one-particle phase space $\mathbb R^{2d}$.
The scope of Lanford's theorem is to prove a law of large numbers for this random field; namely that its distribution concentrates, as $\varepsilon \to 0$, to the deterministic limit driven by the Boltzmann equation:
\begin{equation}
\label{herbert}
\Pi^{\varepsilon}(h,t) \to \int_{\mathbb R^{2d}} f(t) h(x,v)dx dv \;.
\end{equation}

To do this, we are led to study the first and second moments of $\Pi^{\varepsilon}$, that means in particular computing the expectation (with respect to \eqref{eq: initial measure})
\begin{equation}
\label{eq: mean empirical}
{\mathbb E}_\varepsilon \left [ \Pi^{\varepsilon}(h,t) \right] \eqqcolon
\int_{\mathbb R^{2d}} \F_1 (t,x,v) h(x,v)
 \, dxdv
 \end{equation}
 where the right hand side defines the first correlation function $f_1^\varepsilon$.
 Due to the  invariance under permutations of the density distribution $W_N^{\varepsilon}(t)$, it is not hard to see that this amounts to studying its projection on the one-particle space:
\begin{equation}
\label{eq: one particle density}
\F_1  (t,z_1)
= \frac1{\mu_\varepsilon} \,
\sum_{p=0}^{\infty} \,\frac{1}{p!}\, \int dz_{2}\dots dz_{p+1} \,
W_{p+1}^{\varepsilon } (t,\Z_{p+1})  \, .
 \end{equation}
The evolution of $\F_1  (t)$ over time is therefore dictated by that of $W^{\varepsilon }_{N}(t)$ for all integers~$N$.
It  is also easy to see that the equation satisfied by $\F_1$ is not closed :
if we define  the correlation function of order~$k$
  to be the joint probability of $k$ typical, selected particles
 \begin{equation}
\label{eq: mean empirical k}
 \begin{aligned}
&{\mathbb E}_\varepsilon \Big [ \frac1{\mu_\varepsilon ^k} \sum_{(i_1, i_2,\dots i_k)}  h\big(z_{i_1}(t), z_{i_2}(t),\dots z_{i_k}(t)\big)\Big] \\&\eqqcolon
\int_{\mathbb R^{2dk}} \F_k (t, z_1,\dots z_k) h(z_1,\dots , z_k)
 \, dz_1\dots dz_k
 \end{aligned}
 \end{equation}
then   the equation satisfied by $\F_1$ 
   involves the second correlation function $\F_2 $, whose equation itself involves $\F_3 $, and so on. 
We  more  generally obtain that the equation on~$\F_k $ can be written 
\begin{equation}
\label{eq:HSBBGKY}
\partial_t \F_k  +\V_k \cdot \nabla_{\X_k} \F_k  =C^\varepsilon_{k,k+1} \F_{k+1}   \,\,\,\,\,\,\,\,\, \hbox{in } \,\,\,{\mathcal D}^{\varepsilon}_{k}\, ,
\end{equation}
with a specular reflection condition  at the boundary of ${\mathcal D}^{\varepsilon}_{k}$. This set of equations on~$(\F_k)_{k \geq 1}$ is referred to as the BBGKY hierarchy. The collision operator $C^\varepsilon_{k,k+1}$ comes from the boundary term in Green's formula and is defined as follows:
\begin{equation}\label{eq:collopH'}
\begin{aligned}
(C_{k,k+1}^{\varepsilon}    \F_{k+1}) (\Z_k) &\coloneqq \sum_{ i = 1}^{k}   \int    \F_{k+1}   (\Z_k^{\langle i \rangle},x_i,v'_i,x_i+\varepsilon \omega,w') \big( (
w- v_i
)\cdot \omega \big)_+ \, d \omega dw  \\
& \qquad \qquad -   \int    \F_{k+1} (\Z_k^{\langle i \rangle},x_i,v_i,x_i+\varepsilon \omega,w ) \big( (
w - v_i
)\cdot \omega \big)_- \, d \omega dw \, ,
\end{aligned}
\end{equation}
where~$(v'_i,w')$ is recovered from~$(v_i,w)$ through the scattering laws~(\ref{defZ'nij}), and with the notation
$$
\Z_k^{\langle i \rangle} \coloneqq (z_1,\dots,z_{i-1},z_{i+1},\dots,z_k )\,.
$$
Note that if $C_{1,2}$ denotes the formal limit of $C_{1,2}^\varepsilon$ as~$\varepsilon$ goes to zero, then the right-hand side~$Q(f,f)$ of  the Boltzmann equation~(\ref{eq:BCO}) is none other than $C_{1,2}(f^{\otimes 2})$.

The key result   derived by   \textcite{Lanford} is the convergence  for short times of~$\F_1 (t)$
 to~$ f(t)$,   the solution of the Boltzmann equation with initial data $f^0$.

\begin{theo}[\cite{Lanford}]
\label{thm: Lanford}
 Under   assumption~\eqref{eq:weighted space}, there exists a time $T_L >0$ depending on~$C_0$ and~$\beta$ such that, for any  $t \in [0,T_L]$,
\begin{equation}
\label{eq: convergence Boltzmann-Lanford}
\lim_{\varepsilon \to 0} \F_1 (t)  = f(t )  \end{equation}
uniformly in~$\mathbb R^{2d}$.
\end{theo}

\begin{rema}[More on Lanford's convergence result]\label{rmk:more}
 The formula $\eqref{eq: initial measure}$ for the initial data is a specific choice, which is adopted here merely for the sake of simplicity.  \textcite{Lanford,Lanford1} discussed the initial conditions much more thoroughly.
He gave conditions characterising a typical initial state for which Theorem 2.1   holds. 

The convergence~(\ref{eq: convergence Boltzmann-Lanford}) actually holds in a weighted space similar to~(\ref{eq:weighted space}) but with a deteriorated parameter~$\beta$    (indeed~$\F_1 (t)  $ is shown to decay as~$M_{\beta(t)}$ with~$\beta(t)$ decaying linearly in time). 
\end{rema}

\begin{rema}[On the propagation of chaos]
Since the evolution of $\F_1$ involves the whole BBGKY hierarchy, the proof  actually shows that~$\F_k (t)$ converges to~$ f^{\otimes k} (t )$ on the same time interval, almost everywhere in~$\mathbb R^{2dk}$: the initial chaos property \mbox{$\F_k(0) \sim (f^0)^{\otimes k}$} satisfied by  $\eqref{eq: initial measure}$ is therefore asymptotically preserved on~$[0,T_L]$.
 Using that, for \mbox{$k \geq 2$}, the correlation functions are linked to the higher order moments of the field $\Pi^{\varepsilon}$, this implies also the stronger statement $\eqref{herbert}$, i.e.\,the convergence in probability of the empirical measure (see \cite{Spohn91}). 

Nowadays, we are able to control much more detailed information about the random field $\Pi^{\varepsilon}$, such as its local Poissonian statistics, small fluctuations (including a central limit theorem), and large deviations from the typical behavior; see Section~\ref{sct:perspectives} for related considerations.
\end{rema}

\begin{rema}[On the time restriction]
The main drawback of the result is the time restriction~$T_L$: it can be checked that this time is of the same order as the time $T$ on which a fixed point method allows to solve the Boltzmann equation (see Section~\ref{sct:existence}). The   achievement of the main result presented in this paper (Theorem~\ref{thm: DHM} below) is that~$T_L$ can be replaced by the time on which the Boltzmann equation has a smooth solution, regardless of the way such a solution is constructed. See Section~\ref{sct:DHM}   
 for more on this.
\end{rema}

 \subsection{A strategy of the proof}   
 
Lanford's proof, progressively refined and made more quantitative by many authors over the years, is nowadays a classical argument in perturbation theory: we refer for instance to  \textcite{Uchiyama88,Spohn91,CIP94,CGP,Ukai01,GSRT14,  PSS14}. The original proof is built upon the hard-sphere hierarchy \eqref{eq:HSBBGKY}. The guiding idea is: 
\begin{enumerate}
\item[($0'$)] Express the correlation function~$\F_1(t)$ in terms of the initial correlation functions~$\F_k(0)$, using a series expansion in  ``collision histories'', labelled by binary tree graphs, obtained by iterating the Duhamel formula associated with~(\ref{eq:HSBBGKY}) backwards in time. Each term corresponds to the backward flow of a finite group of particles that are successively added in a collisional configuration and then evolve under the reduced hard-sphere dynamics.
  \end{enumerate}
The argument then  proceeds  in three steps.
\begin{enumerate}
\item[($1'$)] Prove that, for sufficiently short times, the series expansion of $\F_1(t)$  constructed in Step ($0'$)  is absolutely convergent, uniformly in $\varepsilon$, in weighted $L^\infty$ spaces such as~(\ref{eq:weighted space}).
\item[($2'$)] Prove that in the low-density scaling \eqref{defmueps}, almost all collision histories have no recollision, meaning that no collision occurs other than the ones prescribed by the collision tree {constructed in Step ($0'$)}. 
\item[($3'$)] The rest of the argument is simple: neglecting recollisions, the  measure of which tends to zero, one sees that the expansion reduces termwise to the iterated Duhamel formula for the Boltzmann equation, and by inspection  that~$\F_k(t) \;\to\; f^{\otimes k}(t)$ when $\varepsilon \to 0$.
\end{enumerate}

The proof we present in Section \ref{sec:proofL} below is not  the original proof by Lanford, though the two proofs share many similarities. The main difference lies in the way the functions~$\F_k(t)$ are explicitly represented. In place of a backward iterated formula, the guiding idea is a {\em cluster expansion} of trajectories as in~\parencite{ICMP}, in turn inspired from previous work on cumulant methods \parencite{PS17,BGSRS23}.  The plan is the following:
\begin{enumerate}
\item[(0)] Implement a {\it dynamical cluster expansion} associated to the hard-sphere dynamics. This leads to express $\F_k(t)$ as an expansion over cluster trajectories labelled by complex 
%\footnote{At first sight, one may believe that replacing Lanford's binary tree graphs by collision graphs is just an innocent combinatorial rearrangement of the same formula. However, there is a major conceptual jump in this step: in the first case, one sums over the number of particles involved, whereas in the second, the sum is taken over the number of collisions (which can be much larger). As we shall see in Section \ref{sct:DHM} the achievement of~\textcite  {DHM} is to control such a formidable sum, over long times, in the low-density limit.})
 {\it collision graphs} (encoding all possible interactions between the particles).
\item[(1)] Obtain uniform bounds on the expansion (as in ($1'$) above).
\item[(2)] Reduce to minimally connected (tree) graphs, by showing that one can discard recollisions. This implies the chaos property as in ($2'$) above.
\item[(3)] Identify the Boltzmann equation after the limit $\varepsilon \to 0$.
\end{enumerate}

Overall, the cluster expansion strategy is more complicated than Lanford's original  method. In fact, it requires to understand and exploit subtle cancellations among graphs with different signs (both in Steps (1) and (3)). The interest of  presenting this strategy is that it prepares for the proof of the main result of this paper in Section \ref{sct:DHM}, which will overcome the time restriction to~$[0,T_L]$ arising from Step (1).

\section{Proof of Lanford's theorem using cluster expansions} 
\label{sec:proofL}

\subsection{A pictorial proof of Steps~(0) to~(3)}
\label{sec: A pictorial proof}

In this paragraph, we provide pictorial explanations of the main steps of the proof of Theorem \ref{thm: Lanford}. For more details, we refer to~\textcite{ICMP}.

\subsubsection*{\textup{\textbf{Step (0)}}}
 Starting from an initial particle configuration $\Z_N (0) = \Z_N$ in $\mathbb{R}^{2dN}$, the collisions in the dynamics induce correlations between the particles. 
Given a time interval $[0,\tau]$, the particle evolution leads to a partition into {(connected) \emph{collision graphs}}~$ \lambda_1, \dots , \lambda_n $ with edges labelled by the particle indices in~$\{1, \dots , N \}$. This is obtained by adding a vertex at the intersection of two particle paths (represented by edges) when they collide in  $[0,\tau]$ {and the $\lambda_i$ represent all the connected components obtained after this construction}; see Figure~\ref{figure: trajectoires [0,tau] no cycle}. 
This representation 
%will be called \emph{collision graph} and
 depends strongly on the time window $[0,\tau]$ and on the initial data. In particular the number~$n$ of   {collision graphs}  ranges from $1$  to $N$.
\begin{figure}[h] 
\centering
\includegraphics[width=5in]{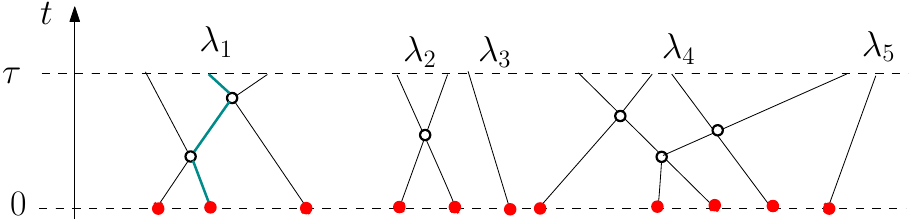} 
\caption{\small This figure encodes the hard-sphere dynamics of 11 particle trajectories during the time interval $[0,\tau]$. The
particles (which are ``independent'') at initial time are represented by red bullets and their paths by lines. A collision between 2 particles is represented by a circle and leads to a scattering of the particles : as an example, a particle trajectory is represented by the green broken line. This evolution leads to the collision graphs $\{ \lambda_1, \dots, \lambda_5 \}$ which contain different numbers of particles (in particular $\lambda_3$ and $\lambda_5$ contain only one particle).
}
\label{figure: trajectoires [0,tau] no cycle}
\end{figure}

Notice that, in general, {collision} graphs are not minimally connected and may include {\it cycles}, as displayed in Figure \ref{figure: trajectoires [0,tau]}. In this case, some of the collisions may be called ``recollisions'' to indicate their role in the generation of the cycles. 

\begin{figure}[h] 
\centering
\includegraphics[width=5in]{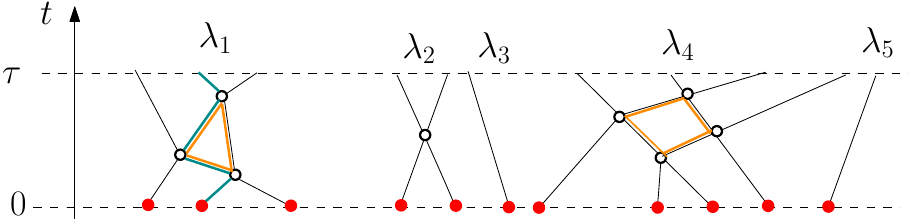} 
\caption{\small Compared  to Figure \ref{figure: trajectoires [0,tau] no cycle}, the collision graphs $\lambda_1$ and $\lambda_4$ display cycles. The top (bottom) collisions in those graphs are recollisions  if we follow the cluster trajectory forward (backward) in time.}
\label{figure: trajectoires [0,tau]}
\end{figure}

Consider now the evolution of a particle which belongs to a collision graph $\lambda_i$. Its trajectory 
 can be complicated due to the collisions, however it is determined only by the particles in $\lambda_i$ (see Figure \ref{figure: trajectoires [0,tau]}). 
 %We will sometimes refer to the particles in $\lambda_i$ as {\em clusters} and to their trajectories  as cluster trajectories.
  The notation $\lambda_i$ will refer not only to a collision graph, but also to the coordinates of the particles indexed by  {the edges of}~$\lambda_i$ and their time evolution.
In this way, there is a one-to-one correspondence between an initial configuration $\{\Z_N \}$
and a decomposition $\{ \lambda_1, \dots , \lambda_n \}$.
This will allow us to change our point of view and,  
 instead of the grand-canonical measure $W^{\varepsilon }_{N}$ with initial datum given in~\eqref{eq: initial measure},
to study a measure on collision graphs 
\begin{equation}
\label{eq: initial measure cluster}
\frac{1}{n!} \widetilde W^{\varepsilon }_{n}(\lambda_1, \dots , \lambda_n) 
\coloneqq \frac{1}{ \widetilde{\mathcal Z}^\varepsilon} \,  \frac{1}{n!}\,\prod_{i=1}^n \nu^\varepsilon ( \lambda_i) \, 
\prod_{i \not = j} {\textbf 1}_{\lambda_i \not \sim \lambda_j} \,,\ \ \  n \geq 0\;,
\end{equation} 
where 
\begin{equation}
\label{nu-def}
 \nu^\varepsilon ( \lambda)\coloneqq {\mu_\eps^{|\lambda|} \over |\lambda|!} (f^0)^{\otimes |\lambda|}  \indc_{\lambda \ \hbox{\tiny is a   collision graph}} \, d\Z_{|\lambda|}\,.\end{equation}
We will not define formally this measure, but intuitively $\widetilde W^{\varepsilon }_{n}$ has a structure similar to $W^{\varepsilon }_{N}$ :
instead of choosing particles with density $f^0$, each collision graph is chosen according to a density $\nu^\varepsilon$ which will determine the number of particles in the graph as well as the corresponding coordinates so that their trajectories satisfy the collision constraints imposed by the graph.
By definition, the dynamics should depend only on the particles within each graph, thus the measure \eqref{eq: initial measure cluster} has also an exclusion condition as the trajectories in different collision graphs $\lambda_i, \lambda_j$ should not {approach at a distance less than~$\eps$}. This is denoted  by $\lambda_i \not \sim \lambda_j$ 
and is the counterpart of the exclusion term~${\textbf 1}_{\mathcal D^{\,\varepsilon}_{\!N}}$ in \eqref{eq: initial measure}.   Finally~$ |\lambda|$ denotes the number of particles in~$\lambda$, which we shall refer to as the {\it size} of~$\lambda$.

In the Boltzmann--Grad scaling~(\ref{defmueps}), a typical particle undergoes one collision per unit time, thus in a small time interval $[0,\tau]$, with $\tau$ smaller than the mean free time between collisions:
\begin{itemize}
\item the size of a collision graph is expected to be of order 1, and therefore the corresponding trajectories should remain simple\footnote{\label{FN2} Here the smallness of $\tau$ is crucial: indeed it has been shown (by numerical and formal arguments) that, for $\tau$ of the order of one mean free time, the size of collision graphs is expected to explode,  by strong analogy with Erd\H{o}s-Renyi random graphs with supercritical parameter \parencite{PSW,PSW2}.}. 
A few collision graphs will be as large as $\log \mu_\varepsilon$, but on the whole, their distribution is well behaved;
\item   the exclusion condition is a weak constraint and will hold outside a  {negligible} set{, meaning that the initial configurations leading to a violation of this condition will have vanishing probability as~$\eps$ goes to zero}. \end{itemize}
The \emph{cluster expansion} method is a well established perturbation theory devised originally to study the 
equilibrium Gibbs distributions of weakly interacting gases. This versatile tool has then been  applied in a variety of cases (see e.g.\,\cite{Poghosyan-Ueltschi}) and we will implement it here to study the weakly interacting collision graphs (which already encode the dynamical information)\footnote{In the context of scaling \eqref{defmueps}, cluster expansion techniques have been studied intensively starting from \textcite{PS17}, and proved to be a natural and powerful strategy when one aims at quantifying correlations and studying fluctuations: see~\textcite{BGSRS20,BGSRS23,BGSRS24,LeBihan25,LeBihan25-1,SW25}.}.

\begin{figure}[h] 
\centering
\includegraphics[width=5in]{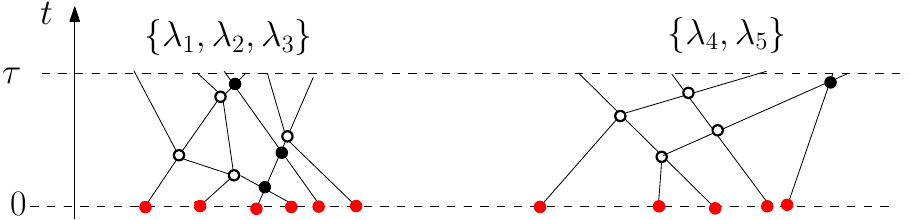} 
\caption{\small 
Above, the collisions graphs in Figure \ref{figure: trajectoires [0,tau]} have been associated to form two clusters : on the left $\lambda_1,\lambda_2, \lambda_3$ are now linked by three overlaps, represented by black dots and on the right $\lambda_4, \lambda_5$ are linked by one overlap. Note that the overlaps on the left generate one cycle.
}
\label{figure: overlaps}
\end{figure}

The  interaction in the measure \eqref{eq: initial measure cluster} comes from the dynamical exclusion condition which we are going to rewrite as 
\begin{equation}
\label{eq: exclusion condition} 
\prod_{i, j \leq n \atop i \not = j} {\textbf 1}_{\lambda_i \not \sim \lambda_j} 
= \prod_{i, j \leq n \atop i \not = j} \big( 1 - {\textbf 1}_{\lambda_i  \sim \lambda_j} \big) \, ,
\end{equation} 
where $\lambda_i  \sim \lambda_j$ means that the collision graphs {\em overlap}, i.e.\;that two particles from each {collision graph} must meet in the time interval $[0,\tau]$, as in Figure \ref{figure: overlaps}. An overlap is not a collision and it does not modify the microscopic dynamics, it is just a convention to represent the complement of  $\{ \lambda_i \not \sim \lambda_j \}$: overlaps represent fictitious dynamics, where particles from different collision graphs can cross each other.
%\begin{figure}[h] 
%\centering
%\includegraphics[width=3in]{pictures/C-O.pdf} 
%\caption{\small The mechanical trajectory of two hard spheres undergoing a collision (left) or an overlap (right). Note that overlaps represent fictitious dynamics.}
%\label{figure: coll vs ov}
%\end{figure}

By expanding the product \eqref{eq: exclusion condition}, we will get terms of the form 
${\textbf 1}_{\lambda_{i_1}  \sim \lambda_{i_2}} {\textbf 1}_{\lambda_{i_2}  \sim \lambda_{i_3}}$ which should be interpreted as follows : the trajectories of the particles are fully determined by each collision graph $\lambda_{i_1} , \lambda_{i_2}, \lambda_{i_3}$ and there are two additional geometric constraints as  $\lambda_{i_1} , \lambda_{i_2}$ and  $\lambda_{i_2}, \lambda_{i_3}$ must overlap
(see Figure \ref{figure: overlaps}).
After this expansion, the relevant objects are the {\em clusters} 
$\sigma_i =  \{ \lambda_{i_1} , \dots , \lambda_{i_\ell} \}$ which are collections of collision graphs constrained to overlap. 
Of course, this operation will possibly lead to further cycles:   this can be seen  in~$\{\lambda_1,\lambda_2, \lambda_3\}$ in Figure \ref{figure: overlaps} for instance, since new cycles have been created due to the overlaps~${\textbf 1}_{\lambda_{1}  \sim \lambda_{2}}  {\textbf 1}_{\lambda_{2}  \sim \lambda_{3}}
{\textbf 1}_{\lambda_{1}  \sim \lambda_{3}} $.

Once the exclusion of the collision graphs has been expanded, there is no more interaction between the clusters. 
Thus the measure \eqref{eq: initial measure cluster} can be upgraded to define a distribution on {\it independent} clusters. This distribution will have a structure reminiscent of a Poisson measure, however it will no longer be a probability measure as the different terms in the expansion of \eqref{eq: exclusion condition}  have signs.
The distribution will have the form
\begin{equation}
\label{eq: initial measure cluster'}
 \frac{1}{ \overline{\mathcal Z}^\varepsilon} \,  \frac{1}{m!}\,\prod_{i=1}^m \tilde \nu^\varepsilon ( \sigma_i) \quad
  m \geq 0\;,
\end{equation} 
where the measure on the cluster $\sigma =  \{ \lambda_1, \dots, \lambda_n\}$  is defined 
  (using \eqref{nu-def}) by
\begin{equation}
\label{tnu-def}
\tilde \nu ^\varepsilon (\sigma)\coloneqq {1\over n!} \varphi_\eps (\lambda_1, \dots, \lambda_n)   \prod_{i=1}^n   \nu^\varepsilon ( \lambda_i) \,,
\end{equation}
and $\varphi_\eps$ encodes the cluster structure between the collision graphs 
$$\varphi_\eps  (\lambda_1, \dots, \lambda_n)  \coloneqq\sum\prod_{\{ j,j'\}  \in E (G_n)} (-\indc_{\lambda_ j \sim \lambda_{j'} })\,,$$
where the sum runs on all connected graphs~$G_n$ with~$n$ vertices, and~$E (G_n)$ denotes the corresponding set of edges. 

\bigskip
Now recall that the Lanford theorem is a statement on~$\F_1(t)$ ---  or more generally on the~$k$-th correlation function~$\F_k(t)$, $k \geq 1$, which describes the behaviour of $k$ typical, selected particles. One thus needs to {\it project} the new measure defined in~(\ref{eq: initial measure cluster'}) by computing correlation functions as in~\eqref{eq: mean empirical k}.

 Let us detail the procedure to recover~$\F_1 (t)$. By \eqref{eq: mean empirical}, using the exchangeability (i.e. the symmetry of the measure under permutations), we have that 
$$
\int_{\mathbb R^{2d}} \F_1 (t,x,v) h(x,v) \, dxdv\approx{\mathbb E}_\varepsilon \left [h(z_1(t)) \right] 
 \,.
$$
The previous decomposition in clusters has two crucial properties
\begin{itemize}
\item The trajectory of a given particle depends only on the particles in the same cluster, in particular $z_1(t)$ is completely prescribed by the cluster $\sigma_1$ containing particle 1;
\item
The measure $W^\varepsilon_N$  is not factorized at the level of particles, but by definition it is at the level of clusters (see \eqref{eq: initial measure cluster'}). 
\end{itemize}
Therefore we can rewrite the previous formula as
$$
\int_{\mathbb R^{2d}} \F_1 (t,x,v) h(x,v)\, dxdv\approx{\mathbb E}_{\tilde\nu^\varepsilon } \left [h(z_1(t)) \right] 
 \,,
$$
keeping only the factor corresponding to the cluster $\sigma_1$, since the average of all other factors in the product (corresponding to clusters independent from particle 1) is 1.
Recall that the  measure $\tilde\nu^\varepsilon (\sigma_1) $ is given by  \eqref{nu-def}-\eqref{tnu-def} encoding the conditioning  by the geometric constraints associated to $\sigma_1$.

\begin{figure}[h] 
\centering
\includegraphics[width=2.8in]{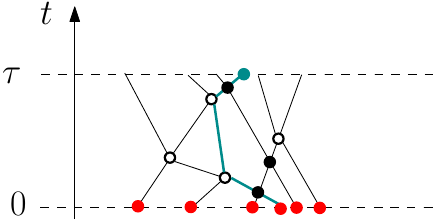} 
\caption{\small Projection leading to the distribution of a typical particle: $\F_1(\tau)$ is expressed as an expansion over cluster trajectories of this form, labelled by {\it connected} collision graphs.
The  particle associated with  $\F_1(\tau)$ is depicted by a green circle at time $\tau$ and its path is also green.
}
\label{figure: equilibrium}
\end{figure}

Figure~\ref{figure: equilibrium} depicts the case of~$\F_1(t)$, where the particle in green is the typical particle, the history of which we are interested in. In this case, only one connected cluster is left: all other clusters disappear from the analysis.
%(after integration in~\eqref{eq: k density function}) since they are independent of the cluster in which this typical particle lies. 

\begin{figure}[h] 
\centering
\includegraphics[width=4in]{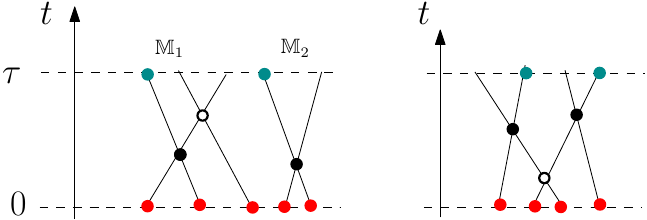} 
\caption{\small Projection leading to the distribution of two typical particles: $\F_2(t)$ is an expansion over cluster trajectories that may belong either to two disconnected collision graphs (on the left), or to the same cluster (on the right).  }
\label{figure: equilibrium k2}
\end{figure}

Of course for~$\F_2(t)$ (more generally~$\F_k(t)$), a similar expansion holds true, but two (or up to $k$) clusters may appear in such an expansion, due to the fact that we select two (or $k$) particles; see Figure~\ref{figure: equilibrium k2}.

From now on and in agreement with Figure~\ref{figure: equilibrium} and the following ones,   we shall call ``green particles"   the particles of interest, meaning those  recorded in~$\F_k(t)$ as opposed to those appearing  due to a    collision/overlap.

\subsubsection*{\textup{\textbf{Step (1)}}}
Using the fact that the gas is dilute and the time $\tau$ is small, one can show that typically the clusters contain only few particles, and few collisions and overlaps. Thus the corresponding distribution \eqref{eq: initial measure cluster} can be studied analytically and it is well behaved. 

More precisely, we want to control the norm of the functions $\F_k(t)$, or (which is a very similar task), we want to understand the behaviour of the partition function $\overline{\mathcal Z}^{\varepsilon}$ in \eqref{eq: initial measure cluster'}.
Note that, to compute these functions, we have to sum over all possible collections of clusters and integrate over all the corresponding cluster trajectories. 

First of all, one can prove an estimate on the density of collision graphs $\nu^\varepsilon$ defined by (\ref{nu-def})  (which, recall, contain no overlaps) of the type
\begin{equation}
\label{eq: borne sur la mesure nu}
\sum_\lambda  \int_{ \mathbb R^{2d|\lambda|}} \left| \nu^\varepsilon(\lambda) \right| = O\left(C_0  \mu_\varepsilon \right),
\end{equation}
  provided that  $\tau \ll C_0^{-1}$ where $C_0$ is defined by \eqref{eq:weighted space}.
This is based on geometric estimates, taking into account the {negligible} volume of configurations which should satisfy the~$|\lambda| - 1$ independent collision constraints imposed by the graph $\lambda$ (roughly, each such constraint is associated to  a volume of order $ \mu_\varepsilon^{-1} \tau$).

Secondly, one has to prove an estimate on the density of clusters $\tilde \nu^\varepsilon$ defined by (\ref{tnu-def}). At this point, we observe that the function $\varphi_\eps$ involves all possible connected graphs on $n$ vertices. This grows badly with $n$: however, these graphs  appear with alternating signs in the
expansion of~$\varphi^\eps$ (each overlap carries a $-1$). A powerful argument due to~\textcite{Penrose} allows to cancel all the graphs with cycles, showing that the original expansion can be bounded from above by a simpler expansion running only over minimally connected graphs ($\ell^{\ell-2}$ graphs involving $\ell$ graphs~$ \{ \lambda_{i_1} , \dots , \lambda_{i_\ell} \}$).
Using this and the geometric estimates, one arrives at a uniform bound of the type 
\begin{equation}
\label{eq: borne sur la mesure tnu}
\sum_\sigma  \int \left| \tilde  \nu^\varepsilon(\sigma) \right| = O\left( C_0 \mu_\varepsilon \right)
\end{equation}
  provided that  $\tau \ll C_0^{-1}$.
This eventually  leads to 
\begin{equation}
\label{eq: borne sur Fk temps court}
\F_k(\tau) \leq C^k \, .
\end{equation}
Recalling that each~$\F_k(\tau) $ has been written as a sum over   all possible collections of clusters, this
 uniform bound allows one to reduce the proof of Theorem \ref{thm: Lanford} to the computation of the limit, graph by graph.

\subsubsection*{\textup{\textbf{Step (2)}}}
 Here we want to prove  the asymptotic propagation of chaos, i.e. that the initial independence of particles (up to exclusion) is preserved, up to an error {vanishing with~$\eps$}, at time~$\tau>0$  recalling that this time is much smaller than the typical time between two collisions. This means that going forward in time from time~$0$ to~$\tau$,  when two particles in the cluster of a green particle meet, they must be independent at the instant before: their histories before that encounter should involve different particles altogether.

\begin{figure}[h] 
\centering
\includegraphics[width=2.5in]{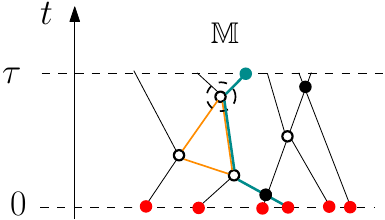} 
\caption{\small A recollision (surrounded by a dashed circle) imposes constraints on the orange part of the trajectory to form a cycle.
}
\label{figure: Fig-recollision-simple}
\end{figure}

This is obviously not the case in  Figure~\ref{figure: Fig-recollision-simple}, since the top collision involves particles that already encountered before. 

\begin{figure}[h] 
\centering
\includegraphics[width=2.8in]{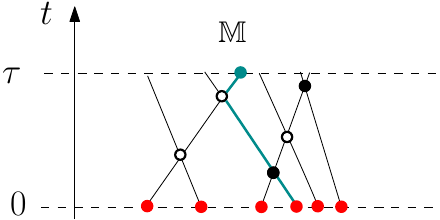} 
\caption{\small Graph (molecule) of leading type in the computation of $f_1^{\varepsilon}(t)$.}
\label{figure: leading cont F1}
\end{figure}

But actually such events are highly improbable. We have already mentioned that at time $0$, each particle in the cluster  is constrained to a {negligible} volume of order $ \mu_\varepsilon^{-1} \tau$. 
Having two particles meet, although they have already  interacted in the past, imposes an even stronger constraint on the velocities and deflection angles at their previous interaction, which leads to an extra decay in $\varepsilon$. One  can prove that in the set of all cluster trajectories of the green particle, the main order  term corresponds to clusters with no cycles at all, meaning trees as in Figure \ref{figure: leading cont F1}. 

%Restricting to such graphs, the limit $ f_1^{\varepsilon}(t) \to f(t)$ for $\varepsilon \to 0$ is easily computed, as $\varepsilon$ is only responsible for small shifts of the trajectories in space, due to the diameter $\varepsilon > 0$ of the hard spheres. 

\bigskip

A very similar argument shows that, in the limit of $f_2^{\varepsilon}(t)$ (see Figure \ref{figure: equilibrium k2} left case), not only are cycles absent, but the two green particles cannot, {asymptotically,} belong to a single connected graph. In other words, the leading contributions come only from clusters with two disconnected components~$(\mathbb{M}_1, \mathbb{M}_2)$, the first containing the green particle 1 and the second containing the green particle 2. As the components $\mathbb{M}_1$ and~$\mathbb{M}_2$ are independent from each other, this leads to the factorisation
$ f_2^{\varepsilon}(t) \to f(t) \otimes f(t) $ for~$\varepsilon \to 0$ (see Step~(3) below).

\begin{rema}
The  technique developed in \parencite{DHM} represents a highly detailed quantitative version of the estimates for clusters with cycles, where it will be shown that, as the number of cycles increases, the corresponding clusters contribute progressively less. 

Consider for instance Figure~\ref{figure: Fig-recollision-simple}, corresponding to a cluster which we denote $\mathbb{M}$. In
the work by Deng, Hani and Ma, clusters of such type are called {\em molecules}. 
They are represented by graphs with edges given by the particles and with vertices associated with collisions or overlaps, with possible cycles generated by collisions (vertices of type $\circ$).
 The orange cycle is the basic example of a \emph{\textup{\{33\}}  molecule} which will be the elementary pattern to gain smallness (we refer to Section~\ref{sec: remainder terms}
for more).
\end{rema}

\subsubsection*{\textup{\textbf{Step (3)}}}
\label{Step3Lanford}  
To conclude, we need to prove that~$\F_1(t)$ is close to the solution of  the Boltzmann equation~(\ref{eq:Beq}), i.e to identify its limit~$f$ as the solution to~(\ref{eq:Beq}).  This means in essence relating Figure \ref{figure: leading cont F1} to the nonlinear PDE~(\ref{eq:Beq}). As mentioned previously, a natural expansion of the solution of~(\ref{eq:Beq}) over trajectories is obtained by simply iterating the Duhamel formula.  Due to the quadratic character of the collision operator, this leads to binary tree graphs. 

\begin{figure}[h] 
\centering
\includegraphics[height=1.5in]{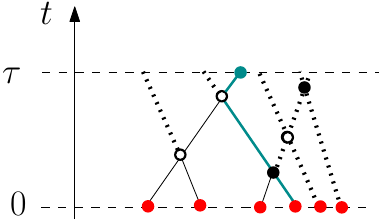} 
\hskip1cm
\includegraphics[height=1.5in]{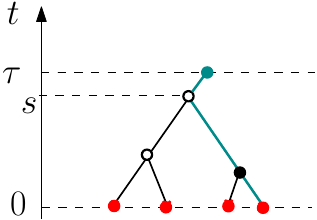} 
\caption{\small On the left, the trajectory of a tagged particle is depicted in green and the part of the cluster which is not involved in this trajectory is represented in dashed lines. 
At the analytical level, the contributions of the dashed lines cancel and the only relevant part of the cluster is the tree represented in the right picture. In the backward direction the tagged trajectory is modified by a collision (at time $s$) determined by two independent trees in $[0,s]$. This is the factorisation property from which Boltzmann equation can be recovered.  
}
\label{eq:cancel B}
\end{figure}
The key point here is to prove that from any tree graph left by Step (2) above (see Figure \ref{figure: leading cont F1}), one can extract a binary tree by following the history of the green particle backward in time, and keeping only the particles that collided directly or indirectly with the green particle, moving backwards. This leads to Figure~\ref{eq:cancel B} where the forgotten pieces of trajectory are represented by dashed lines. 

Is the omission of the dashed part allowed? The answer is affirmative and can be proved by cancellations between graphs with different signs. For instance, in Figure \ref{eq:cancel B} the dashed paths encounter in an overlap $\bullet$, which carries a minus sign. This cancels exactly with the same graph where the $\bullet$ is replaced by a $\circ$ corresponding to a collision.

\subsection{Quantification of the chaos property}\label{sct:chaos}

The chaos property can be made more quantitative by studying more precisely the convergence of $\F_k(t)$ to a tensorised structure. For instance arguing as above, it is not difficult to prove that 
\begin{equation}
\label{eq: propagation of chaos}
\F_2(t) - \F_1 \otimes \F_1(t) = O (\mu_{\varepsilon}^{-1}) 
\end{equation}
for short times, in  $L^1$ norm (see the right part of Figure \ref{figure: equilibrium k2}). A similar estimate is valid 
for 
$
\F_k - (\F_1)^{\otimes k}
$
at least for finite (or slowly diverging) values of $k$.

We would like to introduce now {\em correlation errors} which isolate and quantify the departure from chaos, by means of new functions 
 $E^\varepsilon_k(t)$ with the property that~$E^\varepsilon_k(t)$ vanishes if and only if the $k$ particles are all independent. The decay of $E^\varepsilon_k(t)$ with $\varepsilon$ should control precisely how far the microscopic dynamics remains from the tensorised Boltzmann one, providing a bridge between the existence of the solution $f(t)$ and the validity of the kinetic limit. This method was introduced in \parencite{PS17} and inspired from previous works on stochastic particle systems by \textcite{DMP91, CDMPP91}\footnote{In \parencite{CDMPP91}, a similar method (therein called the method of {\em v-functions}),   proved effective to overcome  Lanford's short time restriction, for a kinetic limit on a lattice.}.

\begin{figure}[h] 
\centering
\includegraphics[height=1.5in]{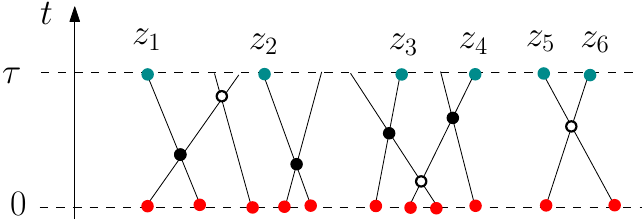} 
\includegraphics[height=1.5in]{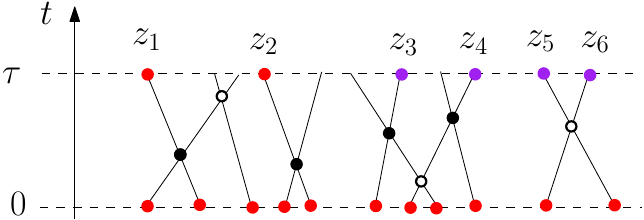} 
\caption{\small In the top figure, representing a contribution to $f^{\varepsilon}_6$, particles~$1,2$ belong to two independent clusters, instead the trajectories of particles~$3,4$ (respectively $5,6$) are correlated by the constraint of belonging to the same cluster. 
Such a configuration cluster will contribute to a term of the type~$\F_1 (t,z_1)\, \F_1 (t,z_2) \, E^\varepsilon_H(t, z_3,z_4,z_5,z_6)$ in \eqref{eq:inverseEJ}
with $H = \{ 3,4,5,6\}$. To indicate this contribution, we update the colours in the bottom figure: red bullets indicate particles which are uncorrelated (their distribution will approximate the solution of the Boltzmann equation), while purple bullets indicate particles associated to the correlation error.
Ultimately the distribution at time $\tau$ is obtained by integrating over the initial variables which are depicted by red bullets and are distributed according to a product measure $f^0$ (as the initial correlations can be ignored by assumption). 
}
\label{eq: E correlations}
\end{figure}

To understand  the correlations intuitively, recall that the $k$-point correlation function~$\F_k(t,z_1, \dots, z_k)$ is determined by projecting on the clusters containing the particles indexed by~$K = \{ 1, \dots, k \}$ at time $t$. At leading order each tagged particle is associated with a different cluster and $\F_k$ factorises  in the low-density limit. 
Correlations emerge between two particles when they belong to the same cluster, see Figure \ref{eq: E correlations}.
The projection on the clusters splits the particles indexed by $K$ into two sets: particles in~$H\subset K$ belong to a cluster containing at least another particle in $H$, whereas each particle in~$K \setminus H$ is associated with a cluster without any other tagged particle. 
The correlation can then be decomposed as 
\begin{equation}
\label{eq:inverseEJ}
\F_k (t, \Z_K) = \sum_{H \subset K}
\bigl(\F_1\bigr)^{\otimes K \setminus H}(t)\,
E^\varepsilon_H (t)\;,
\end{equation}
where $\{ E^\varepsilon_h\}_{h \geq 2}$ are the \emph{correlation errors}. Here we use the abbreviated notation~$\bigl(\F_1\bigr)^{\otimes K \setminus H} (t)= \bigl(\F_1\bigr)^{\otimes (k-|H|)}(t,\Z_{K \setminus H})$, 
and $E^\varepsilon_H (t)= E^\varepsilon_{|H|} (t, \Z_H)$.
The functions~$\F_1$ and~$E^\varepsilon_H$ are obtained by summing over all the clusters compatible with the corresponding dynamical constraints.
Note that the correlation error $E^\eps_{j}$ can be also defined directly from $(\F_\ell)_{\ell \leq j}$ and in particular the first terms are
\[
E^\eps_{2} = \F_2 - \F_1 \!\otimes\! \F_1, \qquad
E^\eps_{3}  = \F_3 - \sum \F_2 \!\otimes\!  \F_1 
+ 2 \F_1 \!\otimes\! \F_1 \!\otimes\! \F_1.
\]
Intuitively, $E^\eps_{j}$ gathers all particles with some correlations, possibly distributed among several clusters: for instance,  in Figure \ref{eq: E correlations}, particles $3,4$ and $5,6$ do not belong to the same cluster and the two pairs are independent.

\bigskip
The correlation errors are therefore much smaller than $
\F_k - (\F_1)^{\otimes k}$ and   \textcite{PS17}     proved that, for suitable constants $\alpha, \gamma \in (0,1)$ and~$C>0$,
\begin{equation}
\label{eq:mainbound}
\int |E^\eps_{H}(t, \Z_K)|\,d\V_H
\le C\, \varepsilon^{\gamma |H|}, 
\qquad  |H|< \varepsilon^{-\alpha},\; 0 \le t < \tau\, .
\end{equation}
Hence, groups of up to $\varepsilon^{-\alpha}$ particles behave independently, up to an exponentially small error.

Once the  first correlation function $\F_1 (\tau)$ has been approximated by the solution of the Boltzmann 
equation $f(\tau)$, one will be able to rewrite \eqref{eq:inverseEJ} up to an error ({vanishing with~$\eps$} in~$L^1$) as 
\begin{equation}
\label{eq:inverseEJ bis}
\F_k (\tau, \Z_K ) = \sum_{H \subset K}
f^{\otimes K \setminus H}(\tau  )\,
E^\varepsilon_H (\tau ) + \text{error} .
\end{equation}
We shall call this a  \emph{cumulant expansion} of~$\F_k(\tau)$. As we shall see in Section \ref{sct:DHM}, a similar cumulant  expansion of the correlation functions provides a convenient ansatz to be propagated to times longer than $\tau$.
Once the tensorised component is known to  approximately satisfy  the Boltzmann equation, and assuming that the solution of the Boltzmann equation is smooth, controlling the growth and decay of the $E^\eps_{j}$'s will be the key.

However, no one knows how to iterate \eqref{eq:mainbound}-\eqref{eq:inverseEJ bis} as a black box, over time steps of length $\tau$. As we shall see  in the next section, Deng, Hani and Ma solve this issue in a drastic way: they complement the ansatz with a {\em full} expansion of the $E-$functions, all the way from time zero up to arbitrary times, bounded by an {\em explicit} graph expansion (see Figure~\ref{fig: time layering}
 and Sections \ref{sec: Strategy of the proof}, \ref{sec: remainder terms}).

 \section{Long time derivation}   \label{sct:DHM}

 \subsection{The theorem by Deng, Hani, Ma}

 The breakthrough of \textcite {DHM} is to extend Lanford's convergence result to much longer times than the mean free time. Note that  this is not a simple technical extension of the previous proof since the series expansions \eqref{eq: borne sur la mesure nu} and \eqref{eq: borne sur la mesure tnu} involving respectively $\nu_\eps $ and $\tilde \nu_\eps$ will be divergent, which is related to the fact that collision graphs on large time intervals $[0,T]$ are expected to have giant components (see footnote \eqref{FN2}).
 The main challenge is then to find a way to combine:
 \begin{itemize}
  \item cluster expansions of the dynamics which converge only locally in time;
 \item the chaos assumption which is known to hold in a strong sense only at time 0. Indeed   \eqref{eq:mainbound} involves a $L^1_v$ norm, in contrast with the assumptions at time zero which guarantee a similar estimate in $L^\infty_v$. The result \eqref{eq:mainbound} cannot be iterated as written, due to this loss in topology\footnote{ It is well known, since the work of Lanford, that the loss of topology is not merely a technical issue but rather an intrinsic feature of the low-density limit, inseparably linked to the transition to irreversibility discussed in Section~\ref{Loschmidt}.}. 
 \end{itemize}
  
 The strategy implemented by Deng, Hani and Ma is therefore to design an iteration scheme which allows to extract at each time step the leading order terms (corresponding to a factorised state) and keep expanding only the remainders. The control on  the leading order terms   
 is then made possible due to a stability argument on the 
 Boltzmann equation, and requires  as explained in Section~\ref{sct:existence} some regularity on its solution, which will be measured in the   weighted space~$L^{\infty,1}_\beta$ defined by the following norm:
 \begin{equation}
\label{eq:weighted space'}
\| f \|_{L^{\infty,1}_\beta} \coloneqq \sum_{k \in \mathbb Z^d} \sup_{|x-k|\leq1} \sup_{v \in   \mathbb R^{d}} \big |f (x,v)\exp( \frac \beta 2 |v|^2)\big | \,.
\end{equation}
 Their result can be stated as follows.
 \begin{theo}[\cite  {DHM}]
\label{thm: DHM}
  In  the low-density scaling \eqref{defmueps},
\begin{equation}
\label{eq: convergence Boltzmann}
\lim_{\varepsilon \to 0} \F_1 (t)  = f(t )  \end{equation}
in~$L^1(\mathbb R^{2d})$, for~$t$ in any finite time interval $[0,T]$ such that the solution $f$ to the Boltzmann equation with initial data $f^0$ satisfies
 $$\forall t \in [0,T] \, , \quad \| f(t) \|_{L^{\infty,1}_\beta}  + \| \nabla_x f (t) \|_{L^{\infty,1}_\beta}  \leq C$$
 for some $C>0$, $\beta>0$.
\end{theo}

Note  that, unlike Theorem~\ref{thm: Lanford},
%(where \eqref{eq: convergence Boltzmann} can be proved to hold uniformly on compact sets) 
 the convergence result~(\ref{eq: convergence Boltzmann})  holds only in $L^1(\mathbb R^{2d})$.

\begin{rema}[On the functional setting]\label{rmk:functional setting'}
The space $L^{\infty,1}_\beta$ encodes both the exponential decay in $v$, and a decay at infinity in $x$. These properties ensure in particular that the total mass is finite, which is not the framework used to derive hydrodynamic limits of the Boltzmann equation (consisting in perturbations of the Maxwellian~$M_\beta$). In a recent preprint \parencite{deng2025-Torus}, the authors 
%announce
{derive} an extension of Theorem~\ref{thm: DHM} to the torus setting.

\end{rema}

\begin{rema}[On diverging times]\label{rmk:diverging time}
Although the convergence result is stated on finite time intervals, the proof shows that it can be extended to slowly diverging times~$T_\eps \to \infty$ as $ \eps \to 0$ if the solution to the Boltzmann equation is global and satisfies
 $$\sup_{t \in \mathbb R^+} \left( \| f(t) \|_{L^{\infty,1}_\beta}  + \| \nabla_x f (t) \|_{L^{\infty,1}_\beta}  \right) <+\infty \,.$$
 Quantifying the  divergence of $T_\eps$  is an important step to study fast relaxation limits leading to hydrodynamic models.
 \end{rema}

\begin{rema}[On the similarities with wave turbulence]
\label{rem: wave}
The proof of Theorem \ref{thm: DHM} shares many similarities with the derivation by the same authors of the wave kinetic equation from the weakly nonlinear Schr\"odinger equation
%As an aside, we should recall that kinetic theory has been very successfully applied to wave equations with small nonlinearities
%and to weakly interacting quantum fluids  \textcite{LS09}.   Deng and Hani started their joint mathematical research with the kinetic limit for the cubic nonlinear Schr\"{o}dinger equation \textcite  {DHw1,DHw2,DH-waves}.
\footnote{{ At the level of weakly nonlinear waves, the second moment of the wave field is governed by the wave kinetic equation, which is the counterpart of the Boltzmann equation. The analogue of the hard sphere empirical density is actually the so-called Wigner function, which is quadratic in the wave fields and random. In the kinetic limit, this object satisfies also a law of large numbers which is the counterpart of \eqref{herbert}.}}  (\cite  {DHw1,DHw2,DH-waves}). Formal derivations of this equation date back to the early work of \textcite  {Peierls} on phonons. A first rigorous result was obtained in \parencite  {LS11} in a close-to-equilibrium setting and for short times, using tools developed by Erd\"os, Salmhofer and Yau in the context of the linear Schr\"odinger equation with random potential (\cite  {EY00,ESY08}). In their breakthrough series of papers, Deng and Hani were able to derive the nonlinear equation and they extended the result to arbitrarily long times. 
We refer to de~\textcite{Suzzoni} for an introduction to these works which we will not comment further on as the setting is very different from the hard sphere dynamics, and we are not familiar enough with the details of the proof.
\end{rema}

\subsection{Strategy of the proof} 
\label{sec: Strategy of the proof}

\begin{figure}[h] 
\centering
\includegraphics[height=1.5in]{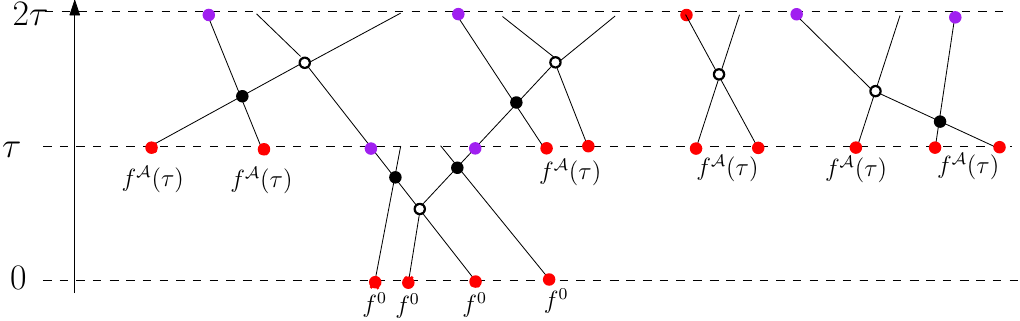} 
\caption{\small  
The collision graph depicted during the time interval $[\tau, 2 \tau]$ involves 
11 particles at time $\tau$. 
It is a contribution to a term of the type~$f^\mathcal A(2\tau) E^\varepsilon_{H}(2\tau)$ with~$|H| = 4$
in an expansion similar to \eqref{eq:inverseEJ bis}. 
We have introduced here the notation $f^\mathcal A$, used by Deng, Hani and Ma to indicate an approximation to the Boltzmann density.
% Its probability can be evaluated by the distribution $\F_{11} (\tau)$ determined in \eqref{eq:inverseEJ bis}. 
 In the example of this figure,~$\F_{11}(\tau)$ is further expanded in a similar way and the purple particles   in $H_1$, with~$|H_1 | =2$, associated to the correlation error at time~$\tau$.
 These particles are dynamically correlated in the past, whereas the particles indicated by red bullets at time~$\tau$ are independently distributed and for those particles we   stop the expansion.
The  error function $E^\eps_{H_1}(\tau)$ is computed by another expansion on $[0,\tau]$ 
with the conditions that the particles in $H_1$ interact dynamically. For simplicity the initial measure is chosen to be a product (neglecting the initial correlations which are much smaller than the dynamical correlations).  }
\label{fig: E correlations 2 fois}
\end{figure}

The strategy devised by \textcite  {DHM} is to establish a cumulant expansion  of the correlation functions  similar to  \eqref{eq:inverseEJ} for further times $T= \mathcal L\tau$, where~$ \mathcal L$  is an integer,
$$\F_k (
\mathcal L\tau,\Z_{K} ) = \sum_{H \subset K}
\bigl(\F_1\bigr)^{\otimes K \setminus H}(\mathcal L\tau )\,
E^\varepsilon_H (\mathcal L\tau )\,,$$
and to provide a graphical representation of the correlation errors $E^\varepsilon_H(\mathcal L \tau)$ by iterating the dynamical cluster expansion \eqref{eq: initial measure cluster'} on the small time intervals $[\ell \tau, (\ell+1) \tau ]$, $\ell < \mathcal L$. The main difference between the dynamical cluster expansion on $[0, \tau]$ and on any further time interval $[\ell \tau, (\ell+1) \tau ]$ is the fact that correlations at time $\ell\tau$ are a priori no longer negligible~: one has to keep track of the dynamical correlations on $[0, \ell \tau]$.
In other words, the representation in Figure \ref{eq: E correlations} where correlations come only from overlaps and collisions during the small time interval has to be changed to take into account previous correlations, meaning that on the 
bottom line some particles will be independent (represented by red bullets
  on the line~$t=\tau$
 in Figure \ref{fig: E correlations 2 fois}), while other will be already correlated (represented by purple bullets in Figure \ref{fig: E correlations 2 fois}).

As already mentioned, iterating the whole cluster expansion on each time step would lead to divergent series (due to the huge number of graphs). Note that this divergence still holds when looking at the iterated expansion of the Boltzmann equation with Duhamel series. {\bf The first  key idea in the work by Deng, Hani and Ma is to extract the leading order terms which are factorised and should compare to the solution of the Boltzmann equation, and to expand only the correlation errors, for which the combinatorics of the graphs should be compensated by the smallness coming from the geometric constraints.} This idea of extracting the maximally factorised state is not present at all in the proof  of Theorem~\ref{thm: Lanford}, where propagation of chaos is obtained as a consequence of the convergence of the correlation functions to the Boltzmann hierarchy. It appeared in subsequent works when establishing more quantitative chaos estimates in \parencite{PS17}, or looking at fluctuations in \parencite{BGSRS23}. But, for hard spheres, it had not been implemented before to obtain information on  the maximally factorised state thanks to the solution of the Boltzmann equation, and to improve the convergence.

\medskip
Let us describe  the strategy in a little more   detail. Set~$s \geq 1$. 
To identify the different correlation errors in  $\F_s ((\ell +1) \tau,\Z_s)$, one  argues iteratively. One  first applies the cluster expansion \eqref{eq: initial measure cluster'}  in $[\ell \tau,( \ell +1)\tau]$ and then evaluates it by using the representation  \eqref{eq:inverseEJ}  of~$(\F_{K_{\ell} } (\ell \tau))_{K_{\ell}}$  according to the following rules    (we have added subscripts~$\ell$ to the sets~$H$ and~$K$ to emphasize that the study takes place on the time step~$[\ell   \tau,(\ell +1) \tau$]):
\begin{enumerate}
\item the variables in $K_{\ell} \setminus H_{\ell}$ are integrated with respect to the density $\F_1 (\ell \tau)$ which is close to $f(\ell\tau)$  by the induction assumption,
\item the variables in $H_{\ell}$ are called \emph{root particles} and they are dynamically correlated in the past history in $[0,\ell\tau]$. Thus another cluster expansion in the previous time interval is needed to compute their joint distribution.
\end{enumerate}

\begin{figure}[h] 
\centering
\includegraphics[height=2in]{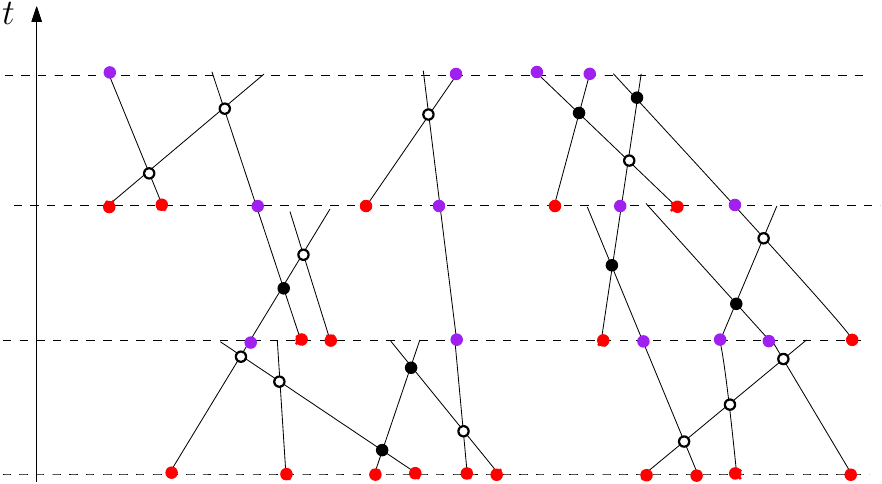} 
\caption{\small  Three time layers representing a correlation term $E_{H} (3 \tau)$  between~$|H| = 4$ root particles at time $3 \tau$.
The dynamical correlations between the root particles form cycles which may spread over several layers.}
\label{fig: time layering}
\end{figure}

The previous simple extraction is actually not sufficient to control the expansion of the correlation errors. Indeed the previous short time cluster expansion  \eqref{eq: initial measure cluster'} leads to collision graphs which can be a priori arbitrarily large and have an unbounded number of recollisions. 
Both events are very unlikely and have been neglected in the time interval $[0,\tau]$ to derive the convergence to the Boltzmann equation (see Section \ref{sec:proofL}). 
For the iteration, it is convenient to remove these events from the main part by introducing {\em truncated microscopic dynamics}\footnote{Section 4.1 in
\textcite{DHM}.}, and deriving  a modified  cumulant expansion of type~\eqref{eq:inverseEJ bis} involving an error term, obtained by a modified Penrose argument\footnote{Section 5 in
\textcite{DHM}.}.

To reach a large time $T= \mathcal L\tau $, the previous strategy is iterated thanks to a \emph{time layering} procedure\footnote{Proposition 3.25 in
\textcite{DHM}.} (see Figure \ref{fig: time layering}).
Applying the rules above at each step, 
one gets the general version of \eqref{eq:inverseEJ bis}
\begin{equation}
\label{eq:inverseEJ general}
\F_k (\ell \tau,\Z_K) = \sum_{H \subset K} (f^\mathcal{A})^{\otimes K \setminus H}(\ell  \tau  ) \, E^\varepsilon_H (\ell  \tau ) + \text{Err}^\eps (\ell  \tau, \Z_K )\,  ,
\end{equation}
where 
\begin{enumerate}
\item $f^\mathcal{A}(\ell\tau)$ is similar to $f(\ell\tau)$ but is a better approximation of $\F_1(\ell\tau)$, and is built on the truncated dynamics. This means that, on each layer, each collision graph has size less than a certain threshold~$\Lambda$ and has less than a fixed number~$\Gamma$ recollisions. Furthermore the number of overlaps in each cluster (or molecule with the terminology of Deng, Hani and Ma) cannot exceed~$\Lambda$. The cut-off $\Lambda \leq | \log \varepsilon |^{C^*_0}$,  for~$C_0^*$ a given constant is allowed to diverge slowly, instead $\Gamma$ will be chosen large but fixed.
\item The correlation error $E^\varepsilon_H$ is obtained by a cluster expansion and also built on the truncated dynamics. 
\item The error~$\text{Err}^\eps $ collects all the remainders accumulating at all steps and will  {converge to zero with~$\eps$} in~$L^1$-norm.
\end{enumerate}

\bigskip

The distribution $f^\mathcal{A} (\ell \tau)$ is approximated step by step by the solution of the Boltzmann equation.  In particular it is shown\footnote{Proposition 6.1, proved in Section 14 in
\parencite{DHM}.} that      
\begin{equation}
\label{eq: approximation Boltzmann-leading}
\| f^\mathcal{A} (\ell \tau ) - f(\ell\tau ) \|_{L^{\infty,1}_{\beta_\ell}} \leq \varepsilon^{\theta_\ell}
\end{equation}
for some~$\beta_\ell>0$  and~$\theta_\ell>0$ that can be explicitly computed\footnote{(1.34) in \parencite{DHM}.}.
The existence of a regular solution to the Boltzmann equation on the time interval $[0,T]$ is  key to establish this estimate, see  Section \ref{sec: The maximally factorised part} below.
Besides the existence of such a solution, the proof of  Deng, Hani and Ma does not use any other mechanism of cancellation in the Boltzmann equation. 

The dynamical correlations are shown to be {negligible} due to the constraint on the re\-collisions\footnote{Proposition 6.2 in
\textcite{DHM}.}:  
\begin{equation}
\label{eq: approximation Boltzmann}
\| E^\varepsilon_H (\ell \tau  )  \|_{L^1} \leq \varepsilon^{3^{-d-2} + c^* |H|},
\end{equation}
where $|H|$ stands for the number of particles in $H$. 
Roughly speaking each root particle  at time~$\mathcal L \tau$ (depicted in purple on Figure \ref{fig: time layering}) is correlated with at least one other root particle and the corresponding dynamical constraints provide a small factor of order~$\varepsilon^\rho$.
As can be seen on Figure \ref{fig: time layering}, this correlation can appear through several layers and many other root particles can be generated in the iterative procedure.
The cut-off~$\Lambda \leq | \log \varepsilon |^{C_0^*}$ on the size of the molecules  (clusters) is needed to tame the combinatorial factor on the number of all molecules (see condition (2) in the description below~\eqref{eq:inverseEJ general}, and Section \ref{sec: remainder terms} below for more). Nevertheless the number of allowed recollisions 
%can be as large as  a power of $|\log \varepsilon |$ leading 
can lead to a huge complexity of dynamical patterns. 
The analysis of the recollisions is the main achievement of  the proof by Deng, Hani and Ma, and is based on a precise algorithmic decomposition of the collision graphs, see Section \ref{sec: remainder terms} below.

 \subsection{The maximally factorised part}
 \label{sec: The maximally factorised part}

 By construction, the distribution $f^\mathcal{A} ((\ell +1) \tau) $ (which is represented by red circles on Figure \ref{fig: E correlations 2 fois} and subsequent figures) is defined recursively by a dynamical cluster expansion on $[\ell \tau, (\ell +1) \tau]$
 \begin{itemize}
 \item  such that each collision graph has size less than $\Lambda$ and has less than $\Gamma$ recollisions, and the number of overlaps in the cluster is less than $\Lambda$,
 \item acting on independent particles distributed according to $f^\mathcal{A} (\ell \tau) $.
 \end{itemize}

On the other hand, the solution  $f((\ell +1) \tau) $ to the Boltzmann equation can be represented by a dynamical cluster expansion on $[\ell \tau, (\ell +1) \tau]$
\begin{itemize}
 \item  such that each cluster is a binary tree graph, 
 \item acting on independent particles distributed according to $f(\ell \tau) $.
 \end{itemize}

 In order to compare both distributions, the idea is to further decompose
\begin{equation}\label{deltaf} 
 f^\mathcal{A}  = \bar f^\mathcal{A}  +\delta \! f^\mathcal{A} 
\end{equation}
 where the leading term $\bar f^\mathcal{A} $ corresponds to the restriction to (bounded) clusters which are minimally connected, and
\begin{equation}\label{deltafA} 
 f = \bar f +\delta \! f 
 \end{equation}
 where the leading term $\bar f$  corresponds to the restriction to (minimally connected) clusters which have size less than $\Lambda$. 
 
 \bigskip
 The difference $\bar f^\mathcal{A} - \bar f$ is controlled using a kind of {\bf weak-strong stability principle for the Boltzmann equation in mild form}.
The cancellation exhibited in Figure~\ref{eq:cancel B} (Step (3)  of the proof of Lanford's theorem described on page~\pageref{Step3Lanford}) shows that the tree graphs involved in the definition of $\bar f^\mathcal{A} ((\ell+1) \tau) $ can be reduced to binary tree graphs.
In other words, $\bar f^\mathcal{A} ((\ell +1)\tau)$ and $\bar f((\ell +1)\tau)$ have the same cluster expansion, acting respectively on $\bar f^\mathcal{A} (\ell \tau)$ and $\bar f(\ell \tau)$, and with a small spatial discrepancy due to the $\eps$ distance separating collisional particles in $\bar f^\mathcal{A} (\ell \tau)$ contrary to~$\bar f  (\ell \tau)$. 

Using the multilinearity of the cluster expansion, we then obtain that each single term in the sum defining $(\bar f^\mathcal{A} - \bar f) ((\ell +1)\tau)$ has at least one occurrence of $(\bar f^\mathcal{A} - \bar f) (\ell \tau)$ or one occurence of $\eps \nabla_x f (\ell \tau) $, which allows to propagate the smallness (up to a small loss due to the contribution of large velocities) from time~$\ell \tau$ to time~$(\ell +1)\tau$.

\bigskip
The contribution  $\delta \! f$ in~(\ref{deltafA}) is  easily shown to be small, as the remainder of a geometric series (recall that the time-step is~$\tau$ only, so the series is under control). The contribution~$\delta \! f^\mathcal{A} $ of the subleading terms in  $f^\mathcal{A} $ in~(\ref{deltaf}) requires much more care because of the need to compensate the loss due to large velocities. This loss is more or less proportional to the size of the molecule, i.e. to the total number of collisions and overlaps which can be much bigger than the number of particles in the cluster. One thus   needs to gain additional smallness, which comes from the geometric constraints associated to the recollisions. This is done thanks to  a simplified version of the algorithms used to control correlation errors sketched in the section below.

\subsection{The remainder terms}
\label{sec: remainder terms}

  The core of the paper by \textcite{DHM} is the estimate \eqref{eq: approximation Boltzmann} of the remainder term which takes into account the recollisions.
Recall that the time interval~$[0,T]$ is split into $\mathcal{L}$ layers of size~$\tau$. As explained in Section \ref{sec: Strategy of the proof}, given $\ell \leq \mathcal{L}$,  the remainder term~$E^\varepsilon_H(\ell  \tau , \Z_H)  $ at time~$\ell \tau$  can be estimated in~$L^1$ by computing, for all compatible molecules~$\bbM$, the cost (smallness) of the constraint that the root  particles (recall they are the  particles in~$H$, represented by purple bullets) have to be correlated, as in Figure \ref{fig: time layering}. {A similar combinatorial classification of correlation terms has been previously developed in the context of wave turbulence to derive the wave kinetic equation from the (microscopic) nonlinear Schr\"odinger equation; see \textcite{DH-waves} and references therein.}

The cost of a molecule is an upper bound on the concatenation of the cluster expansion weights associated with $\bbM$ over $\ell$ layers. A molecule $\bbM$ is interpreted as a graph 
such that a collision or an overlap is coded by a node and the edges are particle trajectories between 2 nodes. These nodes have degree 4 (2 incoming and 2 outgoing particles) and if a particle trajectory stops, or is created between two layers, the corresponding node has degree 1 (represented respectively by free ends or by red dots in Figure \ref{fig: time layering}). 

\bigskip
Let us first estimate the combinatorics of molecules. The number of graphs with degree less than 4, with~$m$ nodes and~$\rho$ independent cycles, is bounded by~$C^m m^\rho$. Roughly, there are less than $C^m$ trees and each cycle boils down to adding a link between two nodes. 
Given $m \leq | \log \varepsilon |^{C_0^*}$, the number of molecules is estimated by 
\begin{equation}
\label{eq: combinatoire molecules}
\# \{ \bbM\, , \,  | \bbM| = m, \rho  \    \text{cycles} \} 
\leq C^m m^\rho \leq C^m \, | \log \varepsilon |^{C_0^* \rho} .
\end{equation}
The logarithmic bound on the molecule size was imposed in the construction of each time layer (see condition (2) below \eqref{eq:inverseEJ general}). This is essential to control the combinatorics, and 
the additional error term $\text{Err}^\eps$ in \eqref{eq:inverseEJ general} takes this cut-off into account. 

\medskip

The divergence in the combinatorial factor \eqref{eq: combinatoire molecules} has to be compensated by the cost of the recollisions associated with the molecules involved in $E^\varepsilon_H$. On a single time layer, we have already seen that one single recollision is enough to provide some smallness and to derive the chaos property \eqref{eq: propagation of chaos}.  More generally quantitative bounds on the correlations of $j$ particles were stated in \eqref{eq:mainbound}.
  The bound  provided in \parencite{DHM} on each given molecule is much stronger: not only does it hold on a large time since~$\ell$ is arbitrary,  but  it provides additional smallness for each cycle.  We will concentrate here on the difficult case when the number~$\rho$ of cycles can be much larger than the cardinality of $|H|$. 
{\bf The major breakthrough by Deng, Hani and Ma\footnote{Eq. (9.53) of their paper, up to some terms which have been omitted for simplicity.} is to establish that there is $\gamma>0$ such that   the cost of a given molecule can be estimated by~$ \tau^{|\bbM|/9} \varepsilon^{\gamma \rho} $, meaning that the gain in  smallness on the cost of the molecule is proportional to the number $\rho$ of cycles}. 
The upper bound \eqref{eq: approximation Boltzmann} on~$
\| E^\varepsilon_H (\ell  \tau )  \|_{L^1}$ follows by combining \eqref{eq: combinatoire molecules}
and this upper bound.

\begin{figure}[h] 
\centering
\includegraphics[height=2.2in]{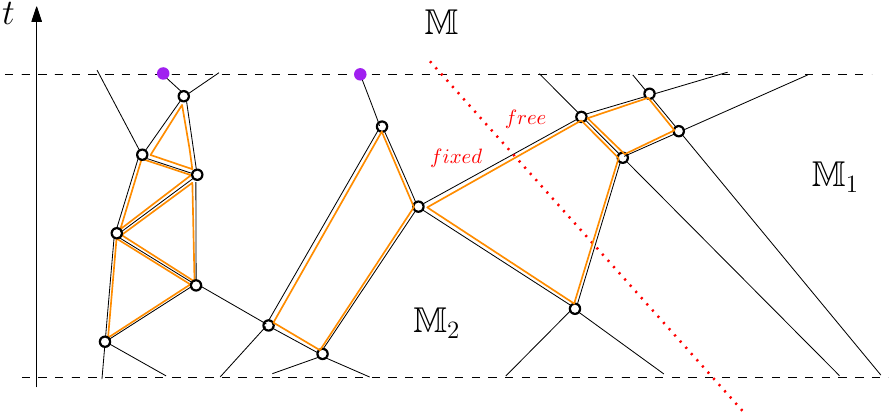} 
\caption{\small  In this figure, the graph represents a molecule $\bbM$ with several cycles and $H$ is of size 2.
On the left a chain of collisions between  a small number of particles  creates several cycles. The  constraints imposed by a molecule $\bbM$ with cycles imply that the parameters of the trajectories take values in a very small set of the phase space. The challenge is to identify the degrees of freedom among these parameters in order to gain a power   $\varepsilon^{\gamma \rho}$ proportional to the number of  cycles $\rho$. 
This is done by splitting the molecule $\bbM$ into smaller entities on which the integration can be performed step by step. For example,  the dotted line decomposes the molecule $\bbM$ into 2 molecules $\bbM_1, \bbM_2$ which will be evaluated one after the other. 
Given a configuration in $\bbM_1$, the particle trajectories crossing the dotted line are frozen  and one has to integrate first over the degrees of freedom restricted to $\bbM_2$ to obtain estimates on the cycles uniformly in $\bbM_1$ and then to integrate over the parameters in $\bbM_1$. 
 A complete splitting procedure is depicted  on Figure \ref{fig: Up algortithm}.}
\label{fig: Chain reaction}
\end{figure}

Given a molecule $\bbM$, computing its cost  
can be seen as  evaluating an integral over the parameters associated with the trajectories compatible with $\bbM$. 
In the proof by Deng, Hani and Ma\footnote{Definition 7.3.}, the molecules are overparametrised as each collision/overlap  is determined by~$8d+1$ variables : one for the position and velocity on each of the 4 edges and one for the time of the collision/overlap. 
Thus a given particle is associated with several sets of distinct variables strongly correlated by the constraints 
at each collision/overlap.
The range of these parameters is further reduced by the dynamical constraints of the cycles imposed by $\bbM$ so that the smallness  follows from the integration over the trajectory parameters
 (see Figure \ref{fig: Chain reaction}).
   The idea behind the duplication of the integration variables is to forget the original structure of the trajectories and to compute the multiple integrals   in any order provided the dynamical constraints are satisfied. 
 For example, the position and velocity of a   particle can be chosen at the beginning of a time layer or in the middle of it. 
 The difficulty is to {\bf determine, among all these variables, which are the degrees of freedom providing the best parametrisation of (a large proportion of) the cycles}\footnote{\label{fn:LS} {This problem has an analogue in the derivation of wave kinetic equations and the strategy is inspired from previous literature on this topic; see \textcite{DH-waves} and references therein.}}. 
  Ordering the parameters will be achieved by an algorithmic procedure described below.
Once an order is determined, one can  use a Fubini argument  to control the cost of the molecule.

\begin{figure}[h] 
\centering
\includegraphics[height=1.2in]{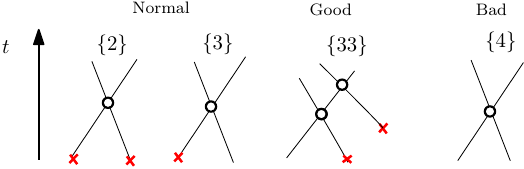} 
\caption{\small 
Some of the elementary molecules are represented above (see Definition 8.11 of \textcite{DHM}).  The crosses represent variables which are given (or fixed).
A $\{2\}$-molecule has 2 incoming fixed edges and therefore the outgoing particle trajectories are fully determined. A $\{3\}$-molecule has only 1 incoming fixed edge and the other incoming edge provides a degree of freedom which can act on the outgoing edges.
A $\{33\}$-molecule encodes the cost of a recollision. Indeed  the single incoming (free) edge (at the bottom) is strongly constrained as the outgoing particle of the first collision has to collide with the 2nd fixed particle trajectory.
Integrating over this constraint leads to a small $\varepsilon^c$  factor and in this sense $\{33\}$ are good molecules. 
Finally, the $\{4\}$-molecule has no constraint and is considered as bad because integrating over too many degrees of freedom leads to a loss.
Similar elementary molecules can be considered by fixing the outgoing parameters instead or by replacing the collisions by overlaps.}
\label{fig: elementary molecules}
\end{figure}

\begin{figure}[h] 
\centering
\includegraphics[height=1.7in]{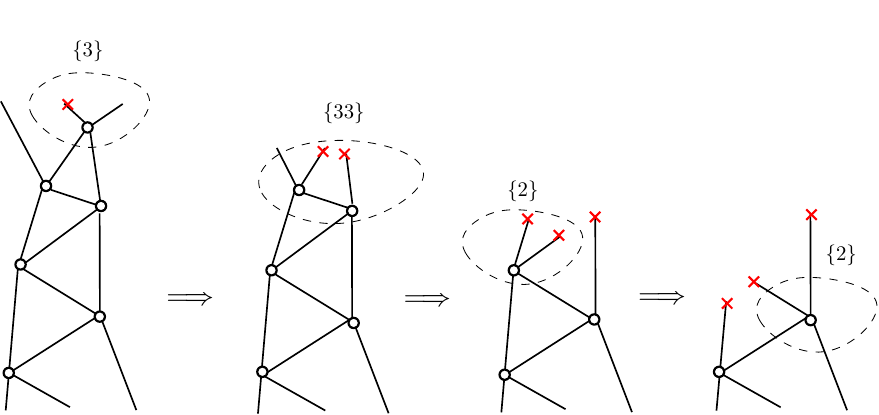} 
\caption{\small  
A series of cycles obtained by multiple recollisions of 3 particles is depicted on the left (see Figure \ref{fig: Chain reaction}). The UP algorithm is applied to the corresponding molecule, starting from the top, in order to decompose it into elementary molecules. After removing the top $\{3\}$-molecule, a $\{33\}$-molecule is extracted in the second step from the first cycle. This exhausts all the degrees of freedom and the next steps will only form $\{2\}$-molecules. Indeed fixing 2 outgoing edges of a collision prescribes exactly the 2 incoming edges  and these constraints propagate along the chain of cycles. Thus only the first constraint can be extracted from this series of recollisions and the gain will not be proportional to the number of cycles.}
\label{fig: Decoupage reaction en chaine}
\end{figure}

As can be imagined from Figure \ref{fig: Chain reaction}, the cycles associated with a molecule can be extremely intricate, leading to long range constraints on the parameters. In particular forcing a cycle to occur may trigger several cycles in  chain reaction, in which case one cannot hope to gain easily a factor $\varepsilon^c$ each time (see also Figure \ref{fig: Decoupage reaction en chaine}). 
 The strategy devised by Deng, Hani and Ma  identifies the degrees of freedom by {\bf decomposing the molecules through several algorithmic procedures\footnote{Sections 10 to 12 in~\parencite{DHM}.}    in order to reduce the complexity of dynamical graphs to elementary molecules}\footnote{Definition 8.11 in \parencite{DHM}.} as in Figure~\ref{fig: elementary molecules} involving only a few parameters.
The elementary molecules are made of one or two collisions or overlaps (called {\it atoms}  by Deng, Hani and Ma). The nature of the elementary molecules is determined by the number of degrees of freedom and they are classified into 3 categories: \emph{normal, good} and \emph{bad}\footnote{See Definition 9.4 of \textcite{DHM} for the precise description of the 11 categories.}. Molecules~$\{3\}$  (and sometimes~$\{2\}$) are the basic elements in the complex molecules and they can be estimated by constants. The dynamical constraints of the cycles are encoded by the $\{33\}$-molecules  which have a cost  of $\varepsilon^c$ (see Figure \ref{fig: elementary molecules}).  
Thus $\{2\}$ and $\{3\}$-molecules are said to be \emph{normal} and~$\{33\}$-molecules are \emph{good}. Molecules $\{4\}$ are \emph{bad} since  all variables are free, and the geometric constraint associated to the collision or overlap is not taken into account.
Note that bad molecules may appear during the decomposition leading to   losses which have to be compensated by enough gains with the $\{33\}$ molecules.
The precise cost of the elementary molecules is evaluated by computing collision integrals depending  only on a few parameters\footnote{Section 9 in \parencite{DHM}.}. 

\bigskip
\bigskip
The  decomposition of a molecule $\bbM$ into elementary molecules $\{2\}, \{3\}, \{33\} $ and~$\{4\}$ is achieved in many steps starting from the large scales 
(see Figure \ref{fig: Chain reaction}) to localise progressively the costly structures of the cycles.  In Figure \ref{fig: Up algortithm} below, an example of a single cycle decomposition is depicted, following the {\bf UP algorithm}\footnote{Definition 12.1 in \parencite{DHM}.}.
\begin{figure}[h] 
\centering
\includegraphics[height=1in]{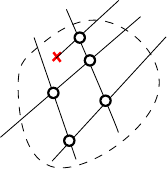} 
$\Rightarrow \atop \text{Step 1}$
\includegraphics[height=1in]{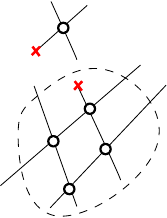} 
$\Rightarrow \atop \text{Step 2}$
\includegraphics[height=1in]{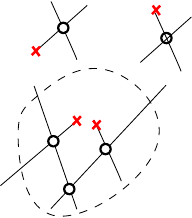} 
$\Rightarrow \atop \text{Step 3}$
\includegraphics[height=1in]{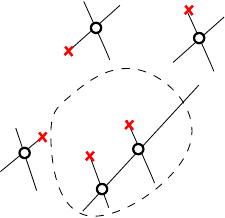} 
\caption{\small  On the left, a molecule involving a cycle is depicted and the edges outside the dashed domain have not been explored yet by the algorithms. In the UP algorithm, the atoms are removed one after the other starting from the top collision in order to keep track of the time ordering associated with the molecule.
In the first step, the top atom  has been extracted as a $\{3\}$-molecule, meaning that 
the remaining edge in the dashed domain is now considered fixed. After the final step, three elementary $\{3\}$-molecules have been identified and a $\{33\}$-molecule will lead to a small cost $\eps^c$.  
 In this example, one  cycle provides exactly one  $\{33\}$-molecule as expected.
For general molecules with overlaps, the decomposition procedure is more delicate as breaking an overlap may lead to non local changes. 
}
\label{fig: Up algortithm}
\end{figure} 
 It is  the simplest algorithm of a series of procedures in order to analyse the numerous structures which can  arise in the molecules.  But it is not enough to  gain a power of $\varepsilon^c$ for each cycle when the cycles in a molecule are very intricate, as in Figure \ref{fig: Decoupage reaction en chaine}.

\newpage
Recall that in \eqref{eq:inverseEJ general}, the main term of the cluster expansion has been devised so that, in  each layer, a cluster has less than $\Gamma$ recollisions. Under this condition, the succession of recollisions between a finite number of particles, as in Figure \ref{fig: Decoupage reaction en chaine}, is no longer a problem as the number of recollisions is bounded by $\Gamma$. Thus extracting a single power~$\varepsilon^c$ is enough to take into account the cycles due to recollisions within a cluster in a given layer, with a cost of~$\eps^{c/\Gamma}$ per cycle. 

\begin{figure}[h] 
\centering
\includegraphics[height=.8in]{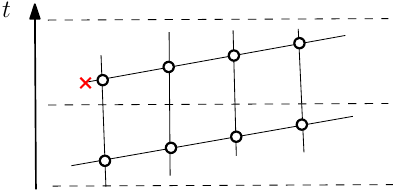} 
$\Rightarrow \atop \text{Step 1}$
\includegraphics[height=.8in]{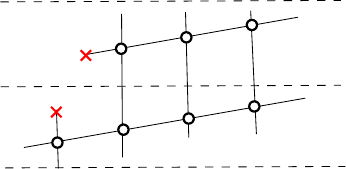} 
$\Rightarrow \atop \text{Step 2}$
\includegraphics[height=.8in]{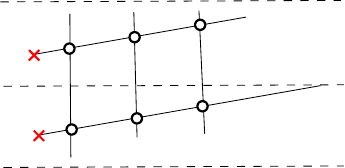} 
\caption{\small  
The left figure depicts a chain of cycles with a ladder structure accross two layers represented by dashed lines.
After removing the first two nodes on the left (as $\{3\}$-molecules), one gets a  $\{33\}$-molecule for each step of the ladder.
This leads to a gain  $\varepsilon^{c\rho} $ proportional to the length of the ladder, i.e. to the number of cycles $\rho$. Note that the order in which the atoms are removed is key and a different procedure could have  led only to a gain  of order $\varepsilon^c$. }
\label{fig: Toy model 1}
\end{figure}

\newpage
However this condition does not prevent the {\bf occurence of an arbitrary number of strongly correlated cycles formed through different layers}.
 There are four main structures of two-layer molecules of this type, named \emph{toy models}\footnote{Section 11 in \parencite{DHM}.} (see Figure~\ref{fig: Toy model 1} for an example of the first toy model).
Each toy model is specific and different algorithms are needed in each case to identify degrees of freedom\footnote{We refer to Section 11 of \textcite{DHM} for an overview of those and to Section 13 for the complete implementation.}.
The toy models have a very streamlined structure~: roughly speaking they look like  two trees on each side of a layer with leaves glued together to form cycles. 

\bigskip
A general molecule $\bbM$ is a priori much more complicated and does not have this tree-like structure. In fact the cycles can spread over several layers and not necessarily over two contiguous layers as in Figure \ref{fig: Toy model 1}.
Given a molecule $\bbM$, an essential part of the proof\footnote{Section 12 in \parencite{DHM}.} is to clean up  $\bbM$ in order to retain only the basic form of a toy model and discard (in a controlled way) the rest of the molecule. This goes through a series of delicate procedures among which the most important are: 
\begin{enumerate}
\item Removing the degeneracies in positions and velocities\footnote{For some particle configurations, a recollision may occur with high probability, e.g. when 3 particles are initially close to each other. The corresponding $\{33\}$-molecules are said to be degenerate and they will not provide the expected gain 
of order $\varepsilon^c$. The occurence of too many such events has to be estimated by other means.};
\item Reducing  the molecule in each layer to a tree by fragmenting each layer;
\item Selecting the two layers across which most of the cycles will form;
\item Cutting exceptional components.
\end{enumerate}
Combining all the previous procedures allows one to select the degrees of freedom associated with a large proportion of the cycles in a molecule.

\begin{rema}
To handle the remainder term in \eqref{eq:inverseEJ general} corresponding to clusters with more than $\Gamma$ recollisions,  a new algorithm had to be devised\footnote{Section 15 of \textcite{DHM}.}. 
The key idea is to use a weak form of dispersion in $\mathbb{R}^d$ provided by the following quantitative statement by \textcite{Burago} : any collection of $q$ hard spheres in $\mathbb{R}^d$ will recollide at most $O( q^{q^2} )$ times. Even though this estimate is far too crude in general, it provides a universal bound, so that if the number of recollisions between $q$ hard spheres is too large, there must be an additional particle involved in the cycles and therefore a new degree of freedom. The use of this dispersion argument is the reason why the   proof does not extend directly to the case of a hard-sphere gas in the torus $\mathbb{T}^d$.
\end{rema}

 \section{Consequences and perspectives}   \label{sct:perspectives}

The paper by \textcite{DHM}  is a crucial advance in deriving macroscopic laws rigorously from microscopic models using the Boltzmann equation as an intermediate description, which was one of the  questions addressed by Hilbert in  his sixth problem. Hydrodynamic models such as the compressible  Euler equations, or the incompressible Navier--Stokes equations can  indeed be recovered as asymptotics of the Boltzmann equation in the fast relaxation limit, i.e. when the collision process is much faster than the transport. In order to complete the program and make a connection with particle systems, one therefore needs to derive the Boltzmann equation for times which are much longer than the mean free time.
 
 As mentioned in Remark \ref{rmk:functional setting'}, although it provides a convergence for times which are much longer than the mean free time, Theorem  \ref{thm: DHM} is not exactly the result we need,  because the functional setting imposed by (\ref{eq:weighted space'}) prevents us from considering initial distributions which are fluctuations around an equilibrium or satisfy at least a bound from below on the macroscopic density, which is required for all hydrodynamic limits. Furthermore, as noted in Remark  \ref{rmk:diverging time}, Theorem  \ref{thm: DHM}  does not provide a quantitative estimate on the time of validity of the Boltzmann approximation, which would be necessary to make explicit the scaling laws leading to the hydrodynamic limits from the particle system via the Boltzmann approximation.
 
 The recent preprint by \textcite{deng2025-Torus} aims at addressing these two objections. It relies on a similar proof, where the use of   \textcite{Burago} is adapted with a new cutting algorithm. It states the convergence from the hard-sphere dynamics to the compressible Euler equations, and incompressible Navier--Stokes equations, under strong regularity assumptions (which are known to be violated in the presence of compressible shocks for instance).
 
 It is worth stressing that the direct hydrodynamic limit from Hamiltonian particle systems, with no smallness in the density, is far from being understood and involves  different kinds of problems. This limit should lead to the Euler  and the Navier--Stokes equations for real fluids\footnote{For a recent discussion on the comparison of the different limit procedures, see e.g.\,\textcite{PS25}.}. At present, a full mathematical result has been obtained only for stochastic dynamics (see~\cite{OVY93,Esposito-Marra-Yau}).

 \bigskip
  A natural question would be to know whether the techniques in \parencite{DHM} can be extended to handle the fluctuations and large deviations both at the kinetic level and at the hydrodynamic level. Using refined cluster expansions (where the correlation error is further decomposed in cumulants corresponding to connected interaction graphs), we have proved in \parencite{BGSRS23} that, in the low-density limit, fluctuations of the hard-sphere dynamics around the Boltzmann approximation can be described by the fluctuating Boltzmann equation, and large deviations can be described by a functional which satisfies some Hamilton-Jacobi equation. However, like Lanford's theorem, these results hold only for short times (except at equilibrium where both the linear fluctuating Boltzmann equation, and the linear fluctuating Stokes-Fourier equations are derived globally in time in \parencite{BGSRS24} and \parencite{BGSRS-AIPH}.
 It has recently been proved by \textcite{Chenjiayue} that the  Hamilton-Jacobi equations are globally well-posed close to equilibrium, so one can hope to push the techniques of \textcite{DHM} to that setting. 

% One important difficulty to extend Thxeorem  \ref{thm: DHM} to the study of fluctuations is that the size of  cumulants decays as powers of $\mu_\eps^{-1}$ , which means that the error terms resulting from the truncation of the dynamics have to be made exponentially small.

 \bigskip
Another important problem, which remains fully open in the field, is the low-density limit of microscopic systems in which the interaction potential is not compactly supported, which is more relevant from the physical point of view
(for some first partial results in this direction, see \cite{Ayi}, and \cite{LeBihan25-1}).
The main difficulty in this case is that the cross section (encoding the rate of the jump process~$(v', v'_\star, \omega) \mapsto (v, v_*, \omega)$ in the Boltzmann collision operator) is not integrable. In particular, the gain and loss terms do not converge separately, and cancellations have to be used to make sense of the equation. 
A natural question is whether this type of cancellations could be exploited indirectly  by assuming the existence of a regular solution, as in the argument by Deng, Hani and Ma presented in Section~\ref{sct:DHM},  to obtain a stable control of the microscopic dynamics. 
%The idea would be to introduce a truncated dynamics with a cut-off on the range of the interactions in the same spirit as in \textcite{Ayi}. The corresponding cross-section would be non uniformly integrable, but maybe with a bound that could be compensated  by the smallness in (\ref{eq: approximation Boltzmann-leading})(\ref{eq: approximation Boltzmann}), provided that the solution of the Boltzmann equation (without cut-off) is smooth and stable.

\printbibliography

@incollection{DolmaireSaffirio,
author = {Dolmaire,T. and Saffirio, C.},
 title = {Derivation of the Boltzmann equation in the disk:quantitative controls of the recollisions},
 booktitle = {Oberwolfach Reports. 20:3},
 pages = {219–2198},
 year = {2024}
 }

@article{KanielShinbrot,
author = {Kaniel, S. and Shinbrot, M.},
		journal = {Communications in Mathematical Physics},
	number = {1},
	pages = {65-84},
	title = {The Boltzmann Equation I. Uniqueness and Local Existence},
	volume = {58},
	year = {1978}}

@incollection{Suzzoni,
 author = {de Suzzoni, A-S.},
 title = {Dérivation de l'équation cinétique associée à l'équation de Schrödinger cubique, d'après Yu Deng et Zaher Hani},
 booktitle = {S\'eminaire Bourbaki. Volume 2023/2024. Expos\'es 1211--1226},
 isbn = {978-2-37905-206-4},
 pages = {341--374},
 year = {2025},
 publisher = {Ast\'{e}risque \textbf{454}},
 language = {French},
 doi = {10.24033/ast.1235},
 zbMATH = {8132295}
}

@article{LeBihan22,
	author = {Le Bihan, C.},
	journal = {Disc. Cont. Dyn Syst.},
	number = {4},
	pages = {1903-1932},
	title = {Boltzmann-Grad limit of a hard sphere system in a box with diffusive boundary conditions},
	volume = {42},
	year = {2022}}

@article{Dolmaire,
	author = {Dolmaire, T.},
	journal = {Kin. Rel. Mod.},
	number = {2},
	pages = {207-268},
	title = {About Lanford's theorem in the half-space with specular reflection},
	volume = {16},
	year = {2023}}

@article{Catapano,
	author = {Catapano, N.},
	journal = {Kin. Rel. Mod.},
	number = {3},
	title = {The rigorous derivation of the Linear Landau equation from a particle system in a weak- coupling limit},
	volume = {11},
	 number = {3},
 pages = {647--695},
	year = {2018}}

@article{GP21,
	author = {Gerasimenko, V. I. and Gapyak, I.V.},
	journal = {Rev. Math. Phys.},
	number = {2},
	title = {Boltzmann-Grad asymptotic behavior of collisional dynamics},
	volume = {33},
	pages = {32 p.},
 note = {Article N. 2130001},
	year = {2021}}

@article{MT12,
	author = {Matthies, K. and Theil, F.},
	journal = {SIAM J. Math. Anal.},
	number = {6},
	pages = {4345-4379},
	title = {A semigroup approach to the justification of kinetic theory},
	volume = {44},
	year = {2012}}

@article{ESY08,
	author = {Erd\H{o}s, L. and Salmhofer, M. and Yau, H. T.},
	journal = {Acta Math.},
	number = {2},
	pages = {211-277},
	title = {Quantum diffusion of the random Schr\"{o}dinger evolution in the scaling limit I. The non-recollision diagrams.},
	volume = {200},
	year = {2008}}

@article{EY00,
	author = {Erd\H{o}s, L. and Yau, H. T.},
	journal = {Commun. Pure Appl. Math.},
	number = {6},
	pages = {667-735},
	title = {Linear Boltzmann equation as the weak coupling limit of a random Schr\"{o}dinger equation},
	volume = {53},
	year = {2000}}

@article{DHw2,
	author = {Deng, Y. and Hani, Z.},
	journal = {Mem. Amer. Math. Soc.},
	title = {Derivation of the wave kinetic equation: full range of scaling laws.},
	note = {To appear},
	year = {2026}}

@article{DHw1,
	author = {Deng, Y. and Hani, Z.},
	journal = {Invent. Math.},
	title = {Full derivation of the wave kinetic equation},
	 volume = {233},
 number = {2},
 pages = {543--724},
	year = {2023}}

@article{LS11,
	author = {Lukkarinen, J. and Spohn, H.},
	journal = {Invent. Math.},
	pages = {79-188},
	title = {Weakly nonlinear Schr\"{o}dinger equation with random initial data},
	volume = {183},
	year = {2011}}

@article{Peierls,
	author = {Peierls, R.E.},
	journal = {Annalen Physik},
	pages = {1055-1101},
	title = {Zur kinetischen Theorie der W{\"a}rmeleitung in Kristallen},
	volume = {3},
	year = {1929}}

@article{PSW2,
	author = {Patterson, R. and Simonella, S. and Wagner, W.},
	journal = {J. Stat. Phys.},
	number = {1},
	pages = {126-167},
	title = {A Kinetic Equation for the Distribution of Interaction Clusters in Rarefied Gases},
	volume = {169},
	year = {2017}}

@article{PSS14,
	author = {Pulvirenti, M. and Saffirio, C. and Simonella, S.},
	journal = {Rev. Math. Phys.},
	number = {2},
	pages = {1-64},
	title = {On the validity of the Boltzmann equation for short range potentials},
	volume = {26},
	year = {2014}}

@article{Esposito-Marra-Yau,
	author = {Esposito, R. and Marra, R. and Yau, H. T.},
	doi = {10.1007/BF02517896},
	fjournal = {Communications in Mathematical Physics},
	issn = {0010-3616},
	journal = {Commun. Math. Phys.},
	language = {English},
	number = {2},
	pages = {395--455},
	title = {Navier-Stokes equations for stochastic particle systems on the lattice},
	volume = {182},
	year = {1996},
	zbl = {0868.60079},
	zbmath = {1002479}}

@misc{Chenjiayue,
	archiveprefix = {arXiv},
	author = {C. Qi},
	eprint = {2409.02805},
	title = {Global Solution of a Functional Hamilton-Jacobi Equation associated with a Hard Sphere Gas},
	url = {https://arxiv.org/abs/2409.02805},
	year = {2024}}

@article{PS25,
	author = {Pulvirenti, M. and Simonella, S.},
	journal = {Mathematics and Mechanics of Complex Systems},
	number = {2},
	pages = {143-154},
	title = {Some Considerations on the Fluid-Dynamical Limit of Particle Systems},
	volume = {14},
	year = {2026}}

@article{OVY93,
	author = {Olla, S. and Varadhan, S.R.S. and Yau, H.-T.},
	journal = {Comm. Math. Phys.},
	pages = {523-560},
	title = {Hydrodynamical Limit for a Hamiltonian System with Weak Noise},
	volume = {3},
	year = {1993}}

@article{Treves,
	author = {Treves, F.},
	journal = {Trans. Amer. Math. Soc},
	pages = {77--92},
	title = {An abstract nonlinear Cauchy-Kowalewska theorem},
	volume = {150},
	year = {1970}}

@article{Loschmidt,
	author = {Loschmidt, J.},
	journal = {Sitzungsberichte der Kaiserlichen Akademie der Wissenschaften Wien, Math. Naturwiss. Klasse},
	pages = {128--142},
	volume = {73},
	year = {1876}}

@article{Fougeres,
	author = {Foug\`eres, F.},
	doi = {10.1007/s10955-024-03353-1},
	fjournal = {Journal of Statistical Physics},
	issn = {0022-4715,1572-9613},
	journal = {J. Stat. Phys.},
	mrclass = {82C22 (35Q20)},
	mrnumber = {4813272},
	number = {10},
	pages = {Paper No. 136, 16 pp.},
	title = {On the derivation of the linear {B}oltzmann equation from the nonideal {R}ayleigh gas},
	url = {https://doi.org/10.1007/s10955-024-03353-1},
	volume = {191},
	year = {2024}}

@article{Ayi,
	author = {Ayi, N.},
	doi = {10.1007/s00220-016-2821-6},
	fjournal = {Communications in Mathematical Physics},
	issn = {0010-3616,1432-0916},
	journal = {Comm. Math. Phys.},
	mrclass = {82C22 (35Q20 60K35)},
	mrnumber = {3607474},
	mrreviewer = {Govind\ Menon},
	number = {3},
	pages = {1219--1274},
	title = {From {N}ewton's law to the linear {B}oltzmann equation without cut-off},
	url = {https://doi.org/10.1007/s00220-016-2821-6},
	volume = {350},
	year = {2017}}

@article{Nirenberg,
	author = {Nirenberg, L.},
	journal = {Jour. Diff. Geom},
	pages = {561--576},
	title = {An abstract form of the nonlinear Cauchy-Kowalewski theorem},
	volume = {6},
	year = {1972}}

@article{SW25,
	author = {Simonella, S. and Winter, R.},
	journal = {Markov Processes and Related Fields},
	title = {Pointwise decay of cumulants in chaotic states at low density},
	note = {To appear.},
	year = {2026}}

@misc{LeBihan25-1,
	author = {Le Bihan, C.},
	howpublished = {arXiv:2408.03597},
	title = {Long time validity of the linearized Boltzmann uncut-off and the linearized Landau equations from the Newton Law},
	year = {2025}}

@article{PSW,
	author = {Patterson, R. and Simonella, S. and Wagner, W.},
	journal = {Phys. D: Nonlin. Phen.},
	pages = {26-32},
	title = {Kinetic theory of cluster dynamics},
	volume = {335},
	year = {2016}}

@article{Lanford1,
	author = {Lanford, III, O. E.},
	journal = {Ast\'{e}risque},
	pages = {117-137},
	title = {On a derivation of the Boltzmann equation},
	volume = {40},
	year = {1976}}

@article{Spohn97,
	author = {Spohn, H.},
	journal = {Pioneering Ideas for the Physical and Chemical Sciences},
	title = {Loschmidt's Reversibility Argument and the H-Theorem},
	year = {1997}}

@article{Carleman,
	author = {Carleman, T.},
	journal = {Acta Math.},
	pages = {369--424},
	title = {Sur la th{\'e}orie de l'equation int{\'e}grodiff{\'e}rentielle de Boltzmann},
	volume = {60},
	year = {1932}}

@article{vBLLS,
	author = {Beijeren, H. van and Lanford, III, O. E. and Lebowitz, J. and Spohn, H.},
	journal = {Journal Stat. Phys.},
	number = {2},
	title = {Equilibrium Time Correlation Functions in the Low--Density Limit},
	volume = {22},
	 pages = {237--257},
	year = {1980}}

@article{BGSR,
	author = {Bodineau, T. and Gallagher, I. and Saint-Raymond, L.},
	journal = {Inventiones mathematicae},
	title = {The Brownian motion as the limit of a deterministic system of hard-spheres},
	volume = {203},
	 number = {2},
 pages = {493--553},
year = {2016}}

@article{IP,
	author = {Illner, R. and Pulvirenti, M.},
	journal = {Commun. Math. Phys.},
	number = {1},
	pages = {143--146},
	title = {Global validity of the Boltzmann equation for two-and three-dimensional rare gas in vacuum: Erratum and improved result},
	volume = {121},
	year = {1989}}

@misc{DH-waves,
	archiveprefix = {arXiv},
	author = {Deng, Y. and Hani, Z.},
	eprint = {2311.10082},
	title = {Long time justification of wave turbulence theory},
	url = {https://arxiv.org/abs/2311.10082},
	year = {2023}}

@article{Burago,
	author = {Burago, D. and Ferleger, S. and Kononenko, A.},
	journal = {Ann. Math. (2)},
	number = {3},
	pages = {695--708},
	title = {Uniform estimates on the number of collisions in semi- dispersing billiards},
	volume = {147},
	year = {1998}}

@article{BGSRS-AIPH,
	author = {Bodineau, T. and Gallagher, I. and Saint-Raymond, L. and Simonella, S.},
	journal = {Ann. Henri Poincar{\'e}},
	pages = {213-234},
	title = {Dynamics of Dilute Gases at Equilibrium: From the Atomistic Description to Fluctuating Hydrodynamics},
	volume = {25},
	year = {2024}}

@article{BGSRS23-1,
	author = {Bodineau, T. and Gallagher, I. and Saint-Raymond, L. and Simonella, S.},
	journal = {Comm. Pure Appl. Math.},
	number = {12},
	title = {Long-Time Correlations for a Hard-Sphere Gas at Equilibrium},
	volume = {76},
	 number = {12},
 pages = {3852--3911},
	year = {2023}}

@misc{DHM,
	archiveprefix = {arXiv},
	author = {Deng, Y. and Hani, Z. and Ma, X.},
	eprint = {2408.07818v3},
	title = {Long time derivation of the Boltzmann equation from hard sphere dynamics},
	url = {https://arxiv.org/abs/2408.07818v3},
	year = {2024}}

@misc{deng2025-Torus,
	archiveprefix = {arXiv},
	author = {Deng, Y. and Hani, Z. and Ma, X.},
	eprint = {2503.01800},
	title = {Hilbert's sixth problem: derivation of fluid equations via Boltzmann's kinetic theory},
	url = {https://arxiv.org/abs/2503.01800},
	year = {2025}}

@book{Cercignani_book,
	author = {Cercignani, C.},
	isbn = {0-19-857064-3},
	language = {English},
	publisher = {Oxford: Oxford University Press},
	title = {Ludwig {Boltzmann}. {The} man who trusted atoms.},
	year = {1998},
	zbl = {1095.01502},
	zbmath = {5049785}}

@article{Denlinger,
	author = {Denlinger, R.},
	doi = {10.1007/s00205-018-1229-1},
	fjournal = {Archive for Rational Mechanics and Analysis},
	issn = {0003-9527},
	journal = {Arch. Ration. Mech. Anal.},
	language = {English},
	number = {2},
	pages = {885--952},
	title = {The propagation of chaos for a rarefied gas of hard spheres in the whole space},
	volume = {229},
	year = {2018},
	zbl = {1397.35164},
	zbmath = {6910330}}

@article{Poghosyan-Ueltschi,
	author = {Poghosyan, S. and Ueltschi, D.},
	doi = {10.1063/1.3124770},
	fjournal = {Journal of Mathematical Physics},
	issn = {0022-2488},
	journal = {J. Math. Phys.},
	language = {English},
	number = {5},
	pages = {Paper No. 053509, 17 pp.},
	title = {Abstract cluster expansion with applications to statistical mechanical systems},
	url = {wrap.warwick.ac.uk/2195/1/WRAP_Ueltschi_abstract_cluster.pdf},
	volume = {50},
	year = {2009},
	zbl = {1187.82009},
	zbmath = {5716519}}

@article{Boltzmann,
	author = {Boltzmann, L.},
	doi = {10.1080/00411459208203923},
	fjournal = {Sitzungsberichte der Akademie der Wissenschaften},
	issn = {0041-1450,1532-2424},
	journal = {Sitzungsberichte der Akademie der Wissenschaften},
	mrclass = {82C40 (76P05)},
	mrnumber = {1165528},
	mrreviewer = {Reinhard\ Illner},
	number = {3},
	pages = {259--276},
	title = {Weitere Studien uber das Warme gleichgenicht unfer Gasmolakular},
	url = {https://doi.org/10.1080/00411459208203923},
	volume = {21},
	year = {1972}}

@incollection{Lanford,
	author = {Lanford, III, O. E.},
	booktitle = {Dynamical systems, theory and applications ({R}encontres, {B}attelle {R}es. {I}nst., {S}eattle, {W}ash., 1974)},
	mrclass = {82.00 (28A65)},
	mrnumber = {479206},
	mrreviewer = {L.\ A.\ Bunimovich},
	pages = {1--111},
	publisher = {Springer, Berlin-New York},
	series = {Lecture Notes in Phys.},
	title = {Time evolution of large classical systems},
	volume = {38},
	year = {1975}}

@article{Cercignani,
	author = {Cercignani, C.},
	journal = {Transport Theory Statist. Phys.},
	pages = {211--225},
	title = {On the Boltzmann equation for rigid spheres},
	volume = {Vol. 2},
	year = {1972}}

@book{CGP,
	author = {Cercignani, C. and Gerasimenko, V. I. and Petrina, D. I.},
	publisher = {Kluwer Academic Publishers, Netherlands},
	title = {Many-Particle Dynamics and Kinetic Equations},
	year = {1997}}

@article{DiPL,
	author = {DiPerna, R. J. and Lions, P.-L.},
	coden = {ANMAAH},
	doi = {10.2307/1971423},
	fjournal = {Annals of Mathematics. Second Series},
	issn = {0003-486X},
	journal = {Ann. of Math. (2)},
	mrclass = {82A40 (35Q20 45K05 76P05)},
	mrnumber = {1014927 (90k:82045)},
	mrreviewer = {Seiji Ukai},
	number = {2},
	pages = {321--366},
	title = {On the {C}auchy problem for {B}oltzmann equations: global existence and weak stability},
	url = {https://proxy.bu.dauphine.fr:443/http/dx.doi.org/10.2307/1971423},
	volume = {130},
	year = {1989}}

@article{Ukai,
	author = {Ukai, S.},
	doi = {10.3792/pja/1195519027},
	fjournal = {Proceedings of the Japan Academy},
	issn = {0021-4280},
	journal = {Proc. Japan Acad.},
	language = {English},
	pages = {179--184},
	title = {On the existence of global solutions of mixed problem for non-linear {Boltzmann} equation},
	volume = {50},
	year = {1974}}

@article{Ukai01,
	author = {Ukai, S.},
	journal = {Japan J. Indust. Appl. Math.},
	pages = {383--392},
	title = {The Boltzmann-Grad Limit and Cauchy-Kovalevskaya Theorem},
	volume = {18},
	year = {2001}}

@article{alexander,
	author = {Alexander, R.},
	journal = {Commun. in Math. Phys.},
	pages = {217--232},
	title = {Time Evolution for Infinitely Many Hard Spheres},
	volume = {49},
	year = {1976}}

@article{Uchiyama88,
	author = {Uchiyama, K.},
	journal = {Hiroshima Math. J.},
	pages = {245--297},
	title = {Derivation of the Boltzmann equation from particle dynamics},
	volume = {18},
	year = {1988}}

@book{CIP94,
	author = {Cercignani, C. and Illner, R. and Pulvirenti, M.},
	doi = {10.1007/978-1-4419-8524-8},
	isbn = {0-387-94294-7},
	mrclass = {82C40 (76-02 76P05 82-02 82B40)},
	mrnumber = {1307620},
	mrreviewer = {Giuseppe Toscani},
	pages = {viii+347},
	publisher = {Springer-Verlag, New York},
	series = {Applied Mathematical Sciences},
	title = {The mathematical theory of dilute gases},
	url = {https://doi.org/10.1007/978-1-4419-8524-8},
	volume = {106},
	year = {1994}}

@book{Spohn91,
	author = {Spohn, H.},
	publisher = {Springer-Verlag, New York},
	series = {Texts and Monographs in Physics},
	title = {Large Scale Dynamics of Interacting Particles},
	year = {1991}}

@book{GSRT14,
	author = {Gallagher, I. and Saint-Raymond, L. and Texier, B.},
	publisher = {EMS},
	series = {Zurich Lect. in Adv. Math},
	title = {From Newton to Boltzmann: hard spheres and short-range potentials},
	year = {2014}}

@article{ICMP,
	author = {Bodineau, T. and Gallagher, I. and Saint-Raymond, L. and Simonella, S.},
	doi = {10.1063/5.0091199},
	fjournal = {Journal of Mathematical Physics},
	issn = {0022-2488,1089-7658},
	journal = {J. Math. Phys.},
	mrclass = {82C40 (35Q70)},
	mrnumber = {4447103},
	mrreviewer = {Cecil\ Pompiliu\ Gr\"unfeld},
	number = {7},
	pages = {Paper No. 073301, 26 pp.},
	title = {Cluster expansion for a dilute hard sphere gas dynamics},
	url = {https://doi.org/10.1063/5.0091199},
	volume = {63},
	year = {2022}}

@article{PS17,
	author = {Pulvirenti, M. and Simonella, S.},
	doi = {10.1007/s00222-016-0682-4},
	fjournal = {Inventiones Mathematicae},
	issn = {0020-9910,1432-1297},
	journal = {Invent. Math.},
	mrclass = {82C22 (35Q20 35Q82)},
	mrnumber = {3608289},
	mrreviewer = {Giuseppe\ Viglialoro},
	number = {3},
	pages = {1135--1237},
	title = {The {B}oltzmann-{G}rad limit of a hard sphere system: analysis of the correlation error},
	url = {https://doi.org/10.1007/s00222-016-0682-4},
	volume = {207},
	year = {2017}}

@article{BGSRS23,
	author = {Bodineau, T. and Gallagher, I. and Saint-Raymond, L. and Simonella, S.},
	doi = {10.4007/annals.2023.198.3.3},
	fjournal = {Annals of Mathematics. Second Series},
	issn = {0003-486X,1939-8980},
	journal = {Ann. of Math. (2)},
	mrclass = {82C40 (60F10 76P05 82C21)},
	mrnumber = {4660136},
	mrreviewer = {Ofer\ Zeitouni},
	number = {3},
	pages = {1047--1201},
	title = {Statistical dynamics of a hard sphere gas: fluctuating {B}oltzmann equation and large deviations},
	url = {https://doi.org/10.4007/annals.2023.198.3.3},
	volume = {198},
	year = {2023}}

@article{BGSRS24,
	author = {Bodineau, T. and Gallagher, I. and Saint-Raymond, L. and Simonella, S.},
	doi = {10.1214/23-aop1656},
	fjournal = {The Annals of Probability},
	issn = {0091-1798,2168-894X},
	journal = {Ann. Probab.},
	mrclass = {82B40 (35Q20 60H15)},
	mrnumber = {4698029},
	mrreviewer = {Alexander\ Orlov},
	number = {1},
	pages = {217--295},
	title = {Long-time derivation at equilibrium of the fluctuating {B}oltzmann equation},
	url = {https://doi.org/10.1214/23-aop1656},
	volume = {52},
	year = {2024}}

@article{BGSRS20,
	author = {Bodineau, T. and Gallagher, I. and Saint-Raymond, L. and Simonella, S.},
	doi = {10.1007/s10955-020-02549-5},
	fjournal = {Journal of Statistical Physics},
	issn = {0022-4715,1572-9613},
	journal = {J. Stat. Phys.},
	mrclass = {35Q20 (60F10 76P05 82C40)},
	mrnumber = {4131018},
	number = {1-6},
	pages = {873--895},
	title = {Fluctuation theory in the {B}oltzmann-{G}rad limit},
	url = {https://doi.org/10.1007/s10955-020-02549-5},
	volume = {180},
	year = {2020}}

@article{LeBihan25,
	author = {Le Bihan, C.},
	doi = {10.1007/s00205-025-02105-z},
	fjournal = {Archive for Rational Mechanics and Analysis},
	issn = {0003-9527,1432-0673},
	journal = {Arch. Ration. Mech. Anal.},
	mrclass = {82C22 (37J99 70F10)},
	mrnumber = {4917079},
	number = {4},
	pages = {Paper No. 40, 72 pp.},
	title = {Long time validity of the linearized {B}oltzmann equation for hard spheres: a proof without billiard theory},
	url = {https://doi.org/10.1007/s00205-025-02105-z},
	volume = {249},
	year = {2025}}

@article{CDMPP91,
	author = {Caprino, S. and De Masi, A. and Presutti, E. and Pulvirenti, M.},
	journal = {Comm. Math. Phys.},
	number = {3},
	pages = {443--465},
	title = {A derivation of the Broadwell equation},
	volume = {135},
	year = {1991}}

@book{DMP91,
	author = {De Masi, A. and Presutti, E.},
	doi = {10.1007/BFb0086457},
	isbn = {3-540-55004-6},
	mrclass = {60K35 (82C05 82C40)},
	mrnumber = {1175626},
	mrreviewer = {Claudio\ Landim},
	pages = {x+196},
	publisher = {Springer-Verlag, Berlin},
	series = {Lecture Notes in Mathematics},
	title = {Mathematical methods for hydrodynamic limits},
	url = {https://doi.org/10.1007/BFb0086457},
	volume = {1501},
	year = {1991}}

@book{Ruelle,
	author = {Ruelle, D.},
	publisher = {W.A. Benjamin Inc., New York},
	title = {Statistical Mechanics. Rigorous Results},
	year = {1969}}

@incollection{Penrose,
	author = {Penrose, O.},
	booktitle = {Statistical mechanics: foundations and applications},
	publisher = {A. Bak ed., Benjamin, New York},
	title = {Convergence of fugacity expansions for classical systems},
	year = {1967}}

@incollection{Grad,
	author = {Grad, H.},
	booktitle = {Handbuch der Physik},
	pages = {205--294},
	publisher = {Springer-Verlag, Berlin-Gottingen-Heidelberg},
	title = {Principles of the kinetic theory of gases},
	year = {1958}}

\end{document}